\documentclass{article}
\usepackage[utf8]{inputenc}
\usepackage{amssymb}
\usepackage{amsmath}
\usepackage{amsfonts}
\usepackage{amsthm}
\usepackage{latexsym}
\usepackage{graphicx}
\usepackage{fullpage}
\usepackage{color,comment}
\usepackage{cite}
\numberwithin{equation}{section}

\newtheorem{theorem}{Theorem}[section]

\newtheorem{lemma}[theorem]{Lemma}

\title{Traveling wave solutions for  two species   competitive chemotaxis systems}
\author{ T. B. Issa\thanks{tissa2@usfca.edu, The Data Institute, University of San Francisco, San Francisco, CA 94117 U.S.A } \and R. B. Salako\thanks{salako.7@osu.edu, Department of Mathematics, The Ohio State University, Columbus, OH 43210 U.S.A } \and  W. Shen\thanks{wenxish@auburn.edu, Department of Mathematics and Statistics Auburn University, Auburn, AL 36849 U.S.A.
} }
\date{}

\begin{document}

\maketitle

\begin{abstract}
   In this paper, we consider  two species chemotaxis systems with Lotka-Volterra competition reaction terms. Under appropriate conditions on the parameters in such a system, we establish the existence of traveling wave solutions of the system connecting  two spatially homogeneous equilibrium solutions  with wave speed greater than some
   critical number $c^*$. We also show the non-existence of such  traveling waves with speed less than some critical number $c_0^*$,  which is independent of the chemotaxis. Moreover, under suitable hypotheses on the coefficients of the reaction terms, we obtain explicit range for the chemotaxis sensitivity coefficients ensuring $c^*=c_0^*$, which implies that the minimum wave speed exists and  is not affected by the chemoattractant.
\end{abstract}


\medskip 


\section{Introduction}

The current work is concerned with the traveling wave solutions of  two species competition chemotaxis systems of the form
\begin{equation}\label{eq:1}
    \begin{cases}
    u_{1,t}=(u_{1,x}-\chi_1 u_1v_x)_x+u_1(1-u_1-au_2),\quad  x\in\mathbb{R}\cr
    u_{2,t}=(du_{2,x}-\chi_2 u_2v_x)_x+ru_2(1-bu_1-u_2),\quad  x\in\mathbb{R}\cr
    0=v_{xx}-\lambda v+\mu_1u_1+\mu_2u_2,\quad  x\in\mathbb{R},
    \end{cases}
\end{equation}
 where $a,b,d,r,\lambda>0$ and  $\mu_i,\chi_i>0$  ($i=1,2$) are positive constants. In \eqref{eq:1},    $u_i(t,x)$, $i=1,2$ denote the density functions of two mobile species living together in the same habitat and competing for some limited resources available in their environment. These two competing species also produce some chemical substance which affects their reproduction dynamics in the sense that each mobile species have tendency to move toward its higher concentration. The density function of the chemical substance is  denoted by $v(t,x)$ and is being produced at the rates $\mu_i$ by the species $i$, for each $i=1,2$. The positive constant $\chi_i$ measures the sensitivity rate by the species $i\in\{1,2\}$ of the chemical substance. The chemical substance has a self degradation rate given by the positive constant $\lambda$. The positive constants $a$ and $b$ measure the interspecific competition between the mobile species. We assume that the first species diffuses at a rate equal to one while the second species diffuses at a rate $d>0$. The positive constant $r$ is the  intrinsic growth rate of the second mobile species.

 When $\chi_1=\chi_2=0$, the dynamics of the chemotaxis model \eqref{eq:1} is governed by the following classical Lotka-Volterra diffusive competition system,
\begin{equation}\label{eq:1-2}
    \begin{cases}
    u_{1,t}=u_{1,xx}+u_1(1-u_1-au_2),\quad x\in\mathbb{R}\cr
    u_{2,t}=du_{2,xx}+ru_2(1-bu_1-u_2),\quad x\in\mathbb{R}.
    \end{cases}
\end{equation}
The asymptotic dynamics of \eqref{eq:1-2} is  of significant research interests and  considerable results have been established in the literature. It is well known that the large time behavior of solutions to \eqref{eq:1-2} is delicately related to the positive constants $a$ and $b$. For example, consider the kinetic system of ODEs associated with \eqref{eq:1-2}, that is,
\begin{equation}\label{eq:1-3}
    \begin{cases}
    u_{1,t}=u_{1}(1-u_1-au_2)\cr
    u_{2,t}=ru_2(1-bu_1-u_2).
    \end{cases}
    \end{equation}
    It is easily  seen that \eqref{eq:1-3} has a trivial equilibrium solution ${\bf e}_0=(0,0)$, and two semi-trivial equilibrium solutions ${\bf e}_1=(1,0)$ and ${\bf e}_2=(0,1)$. When $0<a,b<1$, or $a,b>1$,  \eqref{eq:1-3} has a positive equilibrium given by ${\bf e}_*=(\frac{1-a}{1-ab},\frac{1-b}{1-ab})$.
     The asymptotic dynamics of \eqref{eq:1-3} depends on the strength of the competition coefficients $a$ and $b$.
      We say that {\it  the competition on the first species $u_1$ (resp. on the second species  $u_2$) is strong} if $a>1$ (resp. $b>1$).
     We say that {\it the competition on the first species $u_1$ (resp. on the second species $u_2$) is weak} if $0<a<1$ (resp. $0<b<1$).
     Based on the magnitudes of $a$ and $b$, the following four important cases arise.

     \medskip
\noindent ${\bf (A1)}$ $0<b<1<a$, which is referred to as the {\it strong-weak competition case}.

\smallskip

     \noindent ${\bf (A2)}$  $0<a<1<b$,  which is referred to as the {\it weak-strong competition case}.

     \smallskip

     \noindent ${\bf (A3)}$
     $0<a,b<1$, which is referred to as the {\it weak-weak competition case}.

     \smallskip

     \noindent ${\bf (A4)}$   $a,b>1$, which is referred to as the {\it strong-strong competition case}

\medskip

 Due to biological applications, we are only interested in non-negative solutions of \eqref{eq:1}, \eqref{eq:1-2}, and
 \eqref{eq:1-3}.  Let $(u(t),v(t))$ be a solution of \eqref{eq:1-3} with $u(0)>0$ and $v(0)>0$. The following results are well known. In the case ${\bf (A1)}$,
       $(u_1(t),u_2(t))\to {\bf e}_2$ as $t\to\infty$  and hence the second population outcompetes the first.
  In the case ${\bf (A2)}$,
      $(u_1(t),u_2(t))\to {\bf e}_1$ as $t\to\infty$ and hence  the first population outcompetes the second.
     In the case ${\bf (A3)}$,  $(u_1(t),u_2(t))\to {\bf e}_*$ as $t\to\infty$ and hence both species coexist for all the time.
     In the case ${\bf (A4)}$,   the limit of $( u_1(t),u_2(t))$ depends on the choice of the initial condition $(u_1(0),u_2(0))$.

\smallskip

Let
$$
C_{\rm unif}^b(\mathbb{R})=\{ u\in C(\mathbb{R})\,|\, u\,\, \text{is uniformly continuous and bounded on}\,\, \mathbb{R}\}
$$
with norm $\|u\|_\infty=\sup_{x\in\mathbb{R}}|u(x)|.$
 The above results also hold for the solutions $(u_1(t,x),u_2(t,x))$ of  \eqref{eq:1-2} with initial functions
 $u_i(0,\cdot)\in C_{\rm unif}^b(\mathbb{R})$, $\inf_{x\in\mathbb{R}}u_i(0,x)>0$ $(i=1,2$).
 For example, in the case ${\bf (A2)}$, if $(u_1(t,x),u_2(t,x))$ is a classical solution  of  \eqref{eq:1-2} with
 $u_i(0,\cdot)\in C_{\rm unif}^b(\mathbb{R})$, $\inf_{x\in\mathbb{R}}u_i(0,x)>0$ $(i=1,2$), then
   $$
   \lim_{t\to\infty} (u_1(t,x),u_2(t,x))={\bf e}_1
   $$
   uniformly in $x\in\mathbb{R}$.

   Consider \eqref{eq:1-2}. It is also
     interesting to know the asymptotic behavior of the solutions $(u(t,x),v(t,x))$  with front like initial functions
     $(u(0,x),v(0,x))$, that is, initial functions connecting two equilibrium solutions of \eqref{eq:1-3}. This is
     strongly related to the so called
traveling wave solutions. A {\it traveling wave solution of \eqref{eq:1-2}} is a classical solution of the form $(u_1(t,x),u_2(t,x))=(U_1(x-ct),U_2(x-ct))$ for some constant $c\in\mathbb{R}$, which is called the
 {\it speed of the traveling wave}. A traveling wave solution is said to {\it connect an equilibrium solution ${\bf e}_+$ of \eqref{eq:1-3} at the right end}  if
     \begin{equation}
         \lim_{x\to\infty}(U_1(x),U_2(x))={\bf e}_+.
     \end{equation}
     It is said to {\it connect an equilibrium solution ${\bf e}_-$ of \eqref{eq:1-3}  at the left end} if
     \begin{equation}
         \lim_{x\to-\infty}(U_1(x),U_2(x))={\bf e}_-.
     \end{equation}
The existence of traveling wave solutions of \eqref{eq:1-2} has been extensively studied (see \cite{Fan1,Gar1,Gar2,Hos1,Kan1,Kan2,Kan3,Le1,Li1,fife}). For example, assume $0<a<1$. It is well known that there is a minimum wave speed $c_{\min}\geq  c_0^*:= 2\sqrt{1-a}$ such that \eqref{eq:1-2} has a monotone traveling wave  solution with speed $c$ connecting the equilibrium  solutions ${\bf e}$  and  ${\bf e}_2$   of \eqref{eq:1-3}
 (at the left end and right end, respectively) if and only if $c\geq c_{\min}$,   where ${\bf e}={\bf e}_1$ in the case ${\bf (A2)}$ and ${\bf e}={\bf e}_*$ in the case
 ${\bf (A3)}$. There is no explicit formula available for $c_{\min}$ in the literature and it is  known (see \cite{Huang1}) that it is possible to have the strict inequality $c_{\min}>c_0^*(=2\sqrt{1-a})$. When $c_{\min}=c_0^*$, it is said that the minimum wave speed is {\it linearly determinate}. The works  \cite{Huang2,Li2,Le1} provide sufficient conditions on the parameters to ensure that $c_{\min}$ is linearly determinate. For example, it is proved in \cite[Theorem 2.1]{Le1} that $c_{\min}=c_0^*$
provided the following {\bf (A5)} holds.

\medskip

\noindent {\bf (A5)} $0<d\le 2$ and $(ab-1)r\le (1-a)(2-d)$.

\medskip

The objective of the present work is to investigate in how far
the traveling wave theory for \eqref{eq:1-2} can be extended to the two species chemotaxis system \eqref{eq:1}.
It is clear  that the space-independent  solutions of the chemotaxis system \eqref{eq:1} are the solutions of the  ODE system \eqref{eq:1-3}.
 Note that case ${\bf (A1)}$ and case ${\bf (A2)}$ can be handled similarly.   In the  following, we focus on  case
   ${\bf (A2)}$ and case ${\bf (A3)}$,  and investigate the existence of traveling wave solutions
 of \eqref{eq:1} connecting the unstable equilibrium ${\bf e}_2$ of \eqref{eq:1-3} at the right end and the stable equilibrium ${\bf e}$ of \eqref{eq:1-3} at the left end, where ${\bf e}={\bf e}_1$ in the case ${\bf (A2)}$  and ${\bf e}={\bf e}_*$ in the case
 ${\bf (A3)}$ (see the following subsection for the definition of traveling wave solutions of \eqref{eq:1}).

\subsection{Definition of traveling wave solutions of \eqref{eq:1}}

 Similarly to \eqref{eq:1-2},  a {\it traveling wave solution of \eqref{eq:1}} is a classical solution of the form $(u_1(t,x),u_2(t,x)$, $v(t,x))=(U_1(x-ct),U_2(x-ct),V(x-ct))$ for some constant $c\in\mathbb{R}$, which is called the {\it speed of the traveling wave}. A traveling wave solution $(u_1(t,x),u_2(t,x),v(t,x))=(U_1(x-ct),U_2(x-ct),V(x-ct))$ of \eqref{eq:1}  is said to {\it connect  ${\bf e}_2$ at the right end} if
     \begin{equation}
         \lim_{x\to\infty}(U_1(x),U_2(x))={\bf e}_2,
     \end{equation}
     and to
     {\it  connect ${\bf e}$  at the left end} if
     \begin{equation}
         \lim_{x\to-\infty}(U_1(x),U_2(x))={\bf e}.
     \end{equation}
     A traveling wave solution $(u_1(t,x),u_2(t,x),v(t,x))=(U_1(x-ct),U_2(x-ct),V(x-ct))$ of \eqref{eq:1}  connecting   ${\bf e}_2$ at the right end is {\it nontrivial} if $U_1(x)>0$ for $x\in\mathbb{R}$.

     As far as the chemotaxis model \eqref{eq:1} is concerned, very little is known about the existence of traveling wave solutions.  There are some recent works on the existence and non-existence of traveling wave solutions and spreading speeds of the single species chemotaxis model. In this regards, we refer the reader to the works in \cite{Nadin1,Sa1,Sa2,Sa3,Sa4,Ham1} and references therein.

     We note that  $(u_1(t,x),u_2(t,x),v(t,x))=(U_1(x-ct),U_2(x-ct), V(x-ct))$ is a traveling wave solution of \eqref{eq:1} connecting ${\bf e_2}$ at the right end and  ${\bf e}$ at the left end  if and only if $(U_1(x),U_2(x),V(x))$ is a  steady state solution of the following system
\begin{equation}\label{eq1}
    \begin{cases}
    u_{1,t}=u_{1,xx}+(c-\chi_1v_x)u_{1,x}+u_1(1-\lambda\chi_1v-(1-\chi_1\mu_1)u_1-(a-\chi_1\mu_2)u_{2})\cr
    u_{2,t}=du_{2,xx}+(c-\chi_2v_x)u_{2,x}+ru_2(1-\frac{\lambda\chi_2}{r}v-(1-\frac{\chi_2\mu_2}{r})u_2-(b-\frac{\chi_2\mu_1}{r})u_{1})\cr
    0=dv_{xx}-\lambda v+\mu_1u_1+\mu_2u_2
    \end{cases}
\end{equation}
complemented with the boundary conditions
\begin{equation}\label{boundary-condtion}
    (u_1(-\infty),u_2(-\infty))={\bf e} \quad\text{and}\quad (u_1(\infty),u_2(\infty))={\bf e_2}.
\end{equation}
 Observe that for any solution $(u_1(t,x),u_2(t,x),v(t,x))$ (respectively steady state $(U_1(x),U_2(x),V(x))$) of \eqref{eq:1} (respectively \eqref{eq1}), the third component $v(t,x)$ (respectively $V(x)$) is uniquely determined by the first two components. Hence for the sake of simplicity in the notations, we write ${\bf u}=(u_1,u_2)$ and ${\bf U}=(U_1,U_2)$ for vectors and denote by ${\bf u}(t,x)={\bf U}^c(x-ct)$ traveling wave solutions of \eqref{eq:1} with speed $c\in\mathbb{R}$.
 In the following, we always assume that $u_i(0,\cdot)\in C_{\rm unif}^b(\mathbb{R})$ with $\inf_{x\in\mathbb{R}}u_i(0,x)\ge 0$
 for $i=1,2$.

\subsection{Standing assumptions and notations}

To state the main results of this paper, we introduce some standing assumptions and notations  in this subsection.

It is well known that the solutions $(u_1(t,x),u_2(t,x))$ of \eqref{eq:1-2} with initials $u_i(0,\cdot)\in C_{\rm unif}^b(\mathbb{R})$
 and $0\le u_i(0,\cdot)$ ($i=1,2$) exist  for all $t\ge 0$ and are bounded.
The first standing assumption is on the global existence and boundedness of classical solutions of \eqref{eq:1}.

\medskip

\noindent {\bf (H1)} $1>\chi_1\mu_1$, $r>\chi_2\mu_2$,  $a\ge\chi_1\mu_2$, $br\ge \chi_2\mu_1$.

\medskip

 Note that, when $\chi_1=0$ and $\chi_2=0$, {\bf (H1)} always holds. Note also that \eqref{eq:1} can be written as
  \begin{equation}\label{eq:1-0}
    \begin{cases}
    u_{1,t}=u_{1,xx}-\chi_1v_x u_{1,x}+u_1(1-\lambda\chi_1v-(1-\chi_1\mu_1)u_1-(a-\chi_1\mu_2)u_{2})\cr
    u_{2,t}=du_{2,xx}-\chi_2v_x u_{2,x}+ru_2(1-\frac{\lambda\chi_2}{r}v-(1-\frac{\chi_2\mu_2}{r})u_2-(b-\frac{\chi_2\mu_1}{r})u_{1})\cr
    0=dv_{xx}-\lambda v+\mu_1u_1+\mu_2u_2.
    \end{cases}
\end{equation}
 The assumption {\bf (H1)}  implies that the solutions $(u_1(t,x),u_2(t,x))$ of \eqref{eq:1} with initials $u_i(0,\cdot)\in C_{\rm unif}^b(\mathbb{R})$ and $0\leq u_i(0,\cdot)$ $(i=1,2)$ exist for all
$t>0$ and, moreover,  if $0\le u_i(0,\cdot)\le M_i$ ($i=1,2$), then  $0\leq u_i(t,\cdot)\le M_i$ $(i=1,2$) for all $t\geq 0$, where
\begin{equation}
\label{M1-M2-eq}
M_1=\frac{1}{1-\chi_1\mu_1} \quad {\rm and} \quad M_2:=\frac{r}{r-\chi_2\mu_2}.
\end{equation}

Throughout the rest of this paper, we assume that {\bf (H1)} holds, and that $M_1, M_2$ are the constants defined in \eqref{M1-M2-eq}, and
\begin{equation}
\label{c-zer0-star-eq}
 c_0^*=2\sqrt{1-a}.
\end{equation}

As mentioned in the above, in the case {\bf (A2)}, that is,
 $0<a<1<b$, ${\bf e}_1$ is a stable equilibrium  solution of \eqref{eq:1-2}.
The second standing assumption is on the stability of ${\bf e}_1$ for \eqref{eq:1}.
\medskip

\noindent {\bf (H2)}  $1>\chi_1\mu_1M_1 +aM_2$ and $b\ge 1$.

\medskip

 Note that {\bf (H2)} implies that $a<1$ and $b\ge 1$, and  when $\chi_1=0$ and $\chi_2=0$, ${\bf (H2)}$ reduces to $a<1$ and $b\ge 1$.
We will prove that {\bf (H2)}  implies that ${\bf e_1}$ is  a stable solution of \eqref{eq:1}  (see Theorem \ref{stability-tm}$(i)$).

It is known that, in the case {\bf (A3)}, that is,   $0<a,b<1$,  ${\bf e}_*$ is a stable equilibrium solution of \eqref{eq:1-2}. The third standing assumption is on the stability of ${\bf e}_*$ for \eqref{eq:1}.

\medskip

 \noindent{\bf (H3)} $1>\chi_1\mu_1M_1 +aM_2$ and $r> brM_1+\chi_2\mu_2M_2$.

\medskip

 Note that {\bf (H3)} implies that $a<1$ and $b<1$,  and when $\chi_1=0$ and $\chi_2=0$,   {\bf (H3)} reduces to $a<1$ and $b<1$.
 We will prove that  {\bf (H3)} implies that ${\bf e}_*$ is a stable solution of \eqref{eq:1} (see Theorem \ref{stability-tm}$(ii)$).

 Consider \eqref{eq:1-2}. In the case {\bf (A2)}, it has traveling wave solutions connecting ${\bf e}_2$ at the right end and
 connecting ${\bf e}_1$ at the left end. In the case {\bf (A3)}, it has traveling wave solutions connecting ${\bf e}_2$ at the right end and
 connecting ${\bf e}_*$ at the left end. The next standing assumption is on the existence of traveling wave solutions of \eqref{eq:1}.

 \medskip

 \noindent {\bf (H4)} $(1-a)>r(M_1a(b-\frac{\chi_2\mu_1}{r})-\frac{1}{M_2})_+ +\chi_2(\mu_1M_1+\mu_2M_2)$.

\medskip

Note that, when $\chi_1=0$ and $\chi_2=0$,  {${\bf (H4)}$ becomes $(1-a)>r(ab-1)_+$.}  We will prove that {\bf (H2)}+{\bf (H4)} (resp. {\bf (H3)}+{\bf (H4)}) implies the existence of traveling wave solutions
of \eqref{eq:1} connecting ${\bf e}_2$ at the right end and
 connecting ${\bf e}_1$ at the left end (resp. connecting ${\bf e}_2$ at the right end and
 connecting ${\bf e}_*$ at the left end) for speed $c$ greater than  some number $c^*(\chi_1,\chi_2)(\ge c_0^*)$ (see Theorem \ref{main-exist-thm}).

 As it is mentioned in the above, when {\bf (A5)} holds, the minimal wave speed $c_{\min}$ of \eqref{eq:1-2} is
 linearly determinate, that is, $c_{\rm min}=c_0^*$. The last standing assumption is on the existence and  linear determinacy of
 the minimal wave speed of \eqref{eq:1}.

\medskip
\smallskip

{
\noindent {\bf (H5)}\quad
$\begin{cases}
(1-a)(1+\chi_1\mu_1M_1)<\lambda\cr
 (1-a)(1-(d-1)_+)\ge r\left(M_1\Big(a+\frac{(1-a)\big(\mu_2M_2+a\mu_1M_1+(1-a)\mu_2/\sqrt{\lambda}\big)}{\lambda-(1-a)(1+\chi_1\mu_1M_1)} \Big)(b-\frac{\chi_2\mu_1}{r})-\frac{1}{M_2}\right)_{+}\cr
\quad \qquad\qquad \qquad\qquad\quad + \chi_2\Big(\frac{\lambda }{\lambda-(1-a)} +\frac{\sqrt{1-a})}{\sqrt{\lambda}}\Big)( \mu_1M_1+\mu_2M_2).
\end{cases}
$

\medskip

 It is clear that {\bf (H5)} implies {\bf (H4)}.  When $\chi_1=\chi_2=0$, the first two equations
 in \eqref{eq:1} are independent of $\lambda$ (hence $\lambda$ can be chosen large enough  such that
$1-a<\lambda$) and  {\bf (H5)} reduces to
\begin{equation}\label{zz-zz-1}
r(ab-1)_+\le (1-a)(1-(d-1)_+),
\end{equation}
which holds trivially if {$d<2$ and $ab\le 1$.
  It will be shown that {\bf (H2)}+ {\bf (H5)} (resp. {\bf (H3)}+\bf (H5)}) implies
  the existence and linear determinacy of the minimal wave speed of \eqref{eq:1} connecting
  ${\bf {e}_2}$ at the right end and ${\bf {e}_1}$ at the left end (resp. connecting
  ${\bf {e}_2}$ at the right end and ${\bf {e}_*}$ at the left end)  (see Theorem \ref{tm-min-wave}).



\medskip

In the following, we introduce some standing notations.
 For every $0<\kappa<\kappa_{\max}:=\min\{\sqrt{1-a},\sqrt{\lambda}\}$, let
 \begin{equation}
 \label{c-kappa-eq}
   c_\kappa=\frac{\kappa^2+1-a}{\kappa}
 \end{equation}
 and
  \begin{equation}
  \label{B-lambda-kappa-eq}
  B_{\lambda,\kappa}=\int_{\mathbb{R}}e^{-\sqrt{\lambda}|z|-\kappa z}dz=\frac{1}{\sqrt{\lambda}-\kappa}+\frac{1}{\sqrt{\lambda}+\kappa}=\frac{2\sqrt{\lambda}}{\lambda-\kappa^2}.
  \end{equation}
It is clear that the maps $(0,\sqrt\lambda)\ni \kappa\mapsto B_{\lambda,\kappa} $ and $(0,\sqrt \lambda)\ni \kappa\mapsto \kappa B_{\lambda,\kappa}  $ are monotone increasing. {We define $\kappa_1^*(\chi_1,\chi_2)$  by
\begin{equation}
\label{kappa-1-star-eq}
\kappa_1^*(\chi_1,\chi_2):=\sup\Big\{ \kappa\in(0,\kappa_{\max})\ :\  \kappa B_{\lambda,\kappa}<\frac{2}{\chi_1\mu_1M_1}\Big \}.
\end{equation}
For every $\kappa\in(0,\kappa_1^*(\chi_1,\chi_2))$ we let $f(\kappa,\chi_1,\chi_2)$ denote the positive solution of the algebraic equation
\begin{equation}\label{new-f-de-1}
    \frac{\kappa B_{\lambda,\kappa}}{2}\Big( \mu_2M_2+\mu_1M_1(a+\chi_1f) \Big)+\frac{\mu_2\kappa^2B_{\lambda,\kappa}}{2\sqrt{\lambda}}=f
\end{equation} and define the function $F$ by
\begin{align}
\label{F-2-eq}
F(\kappa,\chi_1,\chi_2)=  & \, \, r\Big(M_1{(a+\chi_1f(\kappa,\chi_1,\chi_2))}(b-\frac{\chi_2\mu_1}{r})-\frac{1}{M_2}\Big)_{+}\nonumber\\
&\,  +\frac{\chi_2({\lambda} B_{\lambda,\kappa}+2\kappa)( \mu_1M_1+\mu_2M_2)}{2\sqrt{\lambda}}-(1-d)\kappa^2.
\end{align}
Let
   \begin{equation}
   \label{kappa-star-eq}
   \kappa^*_{\chi_1,\chi_2}:=\sup\Big\{\tilde{\kappa}\in(0,\kappa_1^*(\chi_1,\chi_2))\, :\, 1-a\geq F(\kappa,\chi_1,\chi_2),\,\, \forall\, 0<\kappa<\tilde{\kappa} \Big\}
   \end{equation}
   and
   \begin{equation}
   \label{c-star-eq}
   c^*=\frac{(\kappa_{\chi_1,\chi_2}^*)^2+1-a}{\kappa_{\chi_1,\chi_2}^*}.
   \end{equation}
Observe that
$$
F(0,\chi_1,\chi_2)= r\big(M_1a(b-\frac{\chi_2\mu_1}{r})-\frac{1}{M_2}\big)_{+}+ \chi_2(\mu_1M_1+\mu_2M_2).
$$
 It is clear that $\kappa^*_1(\chi_1,\chi_2)>0$ and is well defined.  If hypothesis ${\bf (H4)}$ holds, then  $\kappa^*_{\chi_1,\chi_2}$   and $c^*$ are  well defined as well.
 }

\subsection{Statements of the main results}

In this subsection, we  state the  main results of this paper. The first main result is  on the stability of spatially homogeneous  equilibrium solutions of \eqref{eq:1} and is stated in the following theorem.

\begin{theorem}\label{stability-tm} For given $c\in\mathbb{R}$,
let ${\bf u}(t,x;c)$ be a classical solution of \eqref{eq1} satisfying $\inf_{x\in\mathbb{R}}u_1(0,x;c)>0$. Then the following hold.
\begin{itemize}
    \item[(i)] If  ${\bf (H2)}$  holds, then
    $$
    \lim_{t\to\infty}\|{\bf u}(t,\cdot;c)-{\bf e}_1\|_{\infty}=0.
    $$
    \item[(ii)] If  ${\bf (H3)}$ holds, and $\inf_{x\in\mathbb{R}}u_2(0,x;c)>0$,  then
    $$
    \lim_{t\to\infty}\|{\bf u}(t,\cdot;c)-{\bf e}_*\|_{\infty}=0.
    $$
\end{itemize}

\end{theorem}


\smallskip

We have the following two theorems on the existence and nonexistence  of traveling wave solutions of \eqref{eq:1}.

\begin{theorem}\label{main-exist-thm}
 Suppose that  ${\bf (H4)}$ holds. Then for every $c> c^*$, \eqref{eq:1} has a nontrivial traveling solution ${\bf u}(t,x)={\bf U}^{c}(x-ct)=(U_1^c(x-ct),U_2^c(x-ct))$  with speed $c$  connecting ${\bf e}_2$ at right end and satisfying  that
    \begin{equation}\label{asymp-decay-rate}
        \lim_{x\to\infty}\Big|\frac{U_1^c(x)}{e^{-\kappa x}}-1 \Big|=0 \quad \text{and}\quad \lim_{x\to\infty}\Big| \frac{U_2^c(x)-1}{e^{-\kappa x}}-\frac{(\chi_2\mu_2-rb)}{(1-a)+\frac{r}{M_2}-(d-1)\kappa^2-\frac{\chi_2\mu_2\sqrt{\lambda}}{2}B_{\lambda,\kappa}}\Big|=0,
    \end{equation}
    where $\kappa\in(0,\kappa^*_{\chi_1,\chi_2})$ satisfies $c=c_{\kappa}$. Moreover, the following hold.
\begin{enumerate}
    \item [(i)] If ${\bf (H2)}$ holds, then the traveling wave solution ${\bf u}(t,x)={\bf U}^{c}(x-ct)=(U_1^c(x-ct),U_2^c(x-ct))$
    of \eqref{eq:1} also connects
 ${\bf e}_1$ at the left end.

    \item[(ii)] If ${\bf (H3)}$ holds, then the traveling wave solution ${\bf u}(t,x)={\bf U}^{c}(x-ct)=(U_1^c(x-ct),U_2^c(x-ct))$
    of \eqref{eq:1} also connects
 ${\bf e}_*$ at the left end.
\end{enumerate}
\end{theorem}

\begin{theorem}\label{main-non-exist-thm} For any choice of the positive parameters $\chi_i$ and $\mu_i$, $i=1,2$,
there is no nontrivial traveling wave solution ${\bf u}(t,x)={\bf U}^{c}(x-ct)$ of \eqref{eq:1}  with speed $c< c_0^*$ and connecting ${\bf e_2}$ at the right end.
\end{theorem}

 Observe that  Theorem \ref{main-non-exist-thm}  provides a lower bound $c_0^*(=2\sqrt{1-a})$ for the  speeds of traveling wave solutions of \eqref{eq:1}. This lower bound is independent of the chemotaxis sensitivity coefficients $\chi_1$ and $\chi_2$.
  The following theorem shows that this is the greatest lower bound for the  speeds of traveling wave solutions of \eqref{eq:1}.


\begin{theorem}\label{tm-min-wave}
\begin{itemize}
    \item [(i)] If ${\bf (H2)}$ and ${\bf (H5)}$ hold,  then for every $c> c_0^*$,  there is a traveling wave solution  of \eqref{eq:1} with speed $c$ connecting ${\bf e}_1$ and ${\bf e}_2$.
     If, in addition, {$r>2\chi_2\mu_2$} then there is a traveling wave solution of \eqref{eq:1} with speed $c=c_0^*$ connecting  ${\bf e_1}$ and ${\bf e}_2$.

    \item[(ii)] If ${\bf (H3)}$ and ${\bf (H5)}$ hold,  then  for any  $c\ge c_0^*$, there is a traveling wave solution of \eqref{eq:1} with speed $c$ connecting  ${\bf e_*}$ and ${\bf e}_2$.
\end{itemize}
\end{theorem}

We remark that under the conditions of Theorem \ref{tm-min-wave},  the minimum wave speed of \eqref{eq:1} exists,  is linearly determinate,  and  is not affected by the chemotactic effect. In general, it is not known  whether \eqref{eq:1} has a minimal wave speed, and if so, whether both systems \eqref{eq:1} and \eqref{eq:1-2} have the same minimum wave speed. This question is related to whether the presence of the chemical substance slows down or speeds up the minimum wave speed. Note that Theorem \ref{main-non-exist-thm} shows that the presence of the chemical substance doesn't slow down the minimum wave speed of \eqref{eq:1-2} whenever it is linearly determinate. We plan to continue working on this problem in our future works.
  We see from \eqref{zz-zz-1} that Theorem \ref{tm-min-wave} recovers  \cite[Theorem 2.1]{Le1},  which guarantees that the minimum wave speed of \eqref{eq:1-2} is linearly determinate under  hypothesis \eqref{zz-zz-1}.}

We also remark that there are also several interesting works on the dynamics of  solutions to \eqref{eq:1} when considered on bounded domains, see \cite{Lankeit1,Issa1,Miz1,Neg1,Neg2,Neg3,Tello1,Tello2} and the references therein. For example, the works in \cite{Lankeit1,Miz1,Tello2} studied the stability of the equilibria of \eqref{eq:1} on bounded domains with Neumann boundary conditions, while the works \cite{Issa1,Neg1} considered \eqref{eq:1} on bounded domains with some non-local term in the reaction terms.

\medskip

The rest of the paper is organized as follows. In section \ref{sec-for-stability-tm} we present the proof of Theorem \ref{stability-tm}. In section \ref{super-sol-sec} we construct  some super and sub-solutions to be used in the proof of Theorem \ref{main-exist-thm}. Section \ref{sec-for-t1} is devoted to the proof of Theorem \ref{main-exist-thm}. In section \ref{sec-for-t2} we present the proof of Theorem \ref{main-non-exist-thm}. The proof of Theorem \ref{tm-min-wave} is presented in section \ref{sec-for-min-wave}.

\section{Proof of Theorem \ref{stability-tm}}
\label{sec-for-stability-tm}

In this section, we prove Theorem \ref{stability-tm}. Throughout this section $c\in\mathbb{R}$ is an arbitrary fixed number.  We first prove two lemmas, which are fundamental for the proofs of most of our results in the subsequent sections.

\begin{lemma}\label{lem-01} Suppose that ${\bf (H2)}$ holds and let ${\bf u}(t,x;c)$ be a bounded classical  solution of \eqref{eq1} defined for every $x\in\mathbb{R}$ and $t\in\mathbb{R}$. If $\inf_{x,t\in\mathbb{R}}u_1(t,x;c)>0$, then ${\bf u}(t,x;c)\equiv {\bf e}_1$.
\end{lemma}
\begin{proof}
Let
\begin{equation}\label{1-zzz-z11}
l_1=\inf_{t,x\in\mathbb{R}}{u}_1(t,x;c), \quad L_1=\sup_{t,x\in\mathbb{R}}{u}_1(t,x;c), \quad l_2=\inf_{t,x\in\mathbb{R}}{u}_2(t,x;c) \quad  \text{and} \quad\ L_2=\sup_{t,x\in\mathbb{R}}{u}_2(t,x;c).
\end{equation}
Let $M_1$ and $M_2$ be as in \eqref{M1-M2-eq}.
Observe from the comparison principle for elliptic equations that $\mu_1l_1+\mu_2l_2\leq \lambda v(x;{\bf u}(\cdot,\cdot;c))\leq \mu_1L_1+\mu_2L_2$. Hence
\begin{align*}
     u_{1,t}=&{u}_{1,xx}+(c-\chi_1v_x(\cdot;{\bf u})) {u}_{1,x}+(1-(1-\chi_1\mu_1){u}_1-(a-\chi_1\mu_2){u}_2-\chi_1\lambda v(\cdot;{\bf u})){u}_1 \cr
     \leq& {u}_{1,xx}+(c-\chi_1v_x(\cdot;{\bf u})) {u}_{1,x}+(1-\chi_1\mu_1l_1- (1-\chi_1\mu_1){u}_1){u}_1.
\end{align*}
Then by the comparison principle for parabolic equations, we have
$$
{u}_1(t,x;c)\le M_1(1-\chi_{1}\mu_1l_1)_+,\quad \forall\ t,x\in\mathbb{R}.
$$
This implies that  $1>\chi_1\mu_1l_1$  and that
\begin{equation}\label{1-zzz-z6}
    L_1\leq M_1(1-\chi_1\mu_1 l_1).
\end{equation}

Similarly, observe that
\begin{align*}
     u_{1,t}=&{u}_{1,xx}+(c-\chi_1v_x(\cdot;{\bf u})) {u}_{1,x}+(1-(1-\chi_1\mu_1){u}_1-(a-\chi_1\mu_2){u}_2-\chi_1\lambda v(\cdot;{\bf u})){u}_1 \cr
     \geq& {u}_{1,xx}+(c-\chi_1v_x(\cdot;{\bf u})) {u}_{1,x}+(1-\chi_1\mu_1L_1- aL_2- (1-\chi_1\mu_1){u}_1){u}_1.
\end{align*}
Thus, as in the above, we conclude from the comparison  principle for parabolic equations that
\begin{equation}\label{1-zzz-z7}
    l_1\geq M_1(1-\chi_1\mu_1L_1-aL_2).
\end{equation}

Observe also that
\begin{align*}
     u_{2,t}=&d{u}_{2,xx}+(c-\chi_2v_x(\cdot;{\bf u})) {u}_{2,x}+(r-\frac{r}{M_2}{u}_2-(br-\chi_2\mu_1){u}_1-\chi_2\lambda v(\cdot;{u})){u}_2 \cr
     \leq& d{u}_{2,xx}+(c-\chi_2v_x(\cdot;{\bf u})) {u}_{2,x}+r(1-bl_1- \frac{1}{M_2}{u}_2){u}_2.
\end{align*}
We conclude from the comparison principle for parabolic equations that
\begin{equation}\label{1-zzz-z8}
L_2\leq M_2(1-bl_1)_+.
\end{equation}

From this point, we distinguish two cases and show that $L_1=l_1=1$ and $L_2=0$.

\smallskip

\noindent {\bf Case 1}. $(1-bl_1)_+=0$. In this case,  by \eqref{1-zzz-z8}, $L_2=0$. Thus, taking the difference of \eqref{1-zzz-z6} and \eqref{1-zzz-z7} side-by-side yields
$$
(1-\chi_1\mu_1M_1)(L_1-l_1)\leq 0.
$$
Since $1>\chi_1\mu_1M_1$ (see {\bf (H2)}) and $l_1\leq L_1$, we obtain that $L_1=l_1$, which combined with \eqref{1-zzz-z6} and \eqref{1-zzz-z7} and the fact $L_2=0$ yield $l_1=L_1=1$.

\smallskip

\noindent {\bf Case 2.}  $1-bl_1>0$. In this case,  it follows from \eqref{1-zzz-z6}-\eqref{1-zzz-z8} that
\begin{align*}
l_1\geq& M_1\left(1-\chi_1\mu_1M_1(1-\chi_1\mu_1l_1)-aM_2(1-bl_1) \right)\cr
=& M_1(1-\chi_1\mu_1M_1-aM_2)+M_1((\chi_1\mu_1)^2M_1+abM_2)l_1.
\end{align*}
That is
\begin{equation*}
\left(\frac{1}{M_1}-(\chi_1\mu_1)^2M_1-abM_2 \right)l_1\geq 1-\chi_1\mu_1M_1  -aM_2.
\end{equation*}
Using the fact that $\chi_1\mu_1+\frac{1}{M_1}=1$, a simple computation shows that $
\frac{1}{M_1}-(\chi_1\mu_1)^2M_1=1-\chi_1\mu_1M_1$, and hence
$$ (1-\chi_1\mu_1M_1-abM_2)l_1\ge 1-\chi_1\mu_1M_1-aM_2.$$
This implies that $l_1\geq 1$ since $b\geq 1$ (see {\bf (H2)} holds). Thus, we get from \eqref{1-zzz-z8} that $L_2=0$, so we can proceed as in {\bf case 1} to show that  $l_1=L_1=1$ as well.

From both cases and the definition of $l_1,L_1$, and $L_2$, we obtain ${\bf u}(t,x;c)\equiv (1,0)$, which completes the proof of the lemma.
\end{proof}

\begin{lemma}\label{lem-02}
Suppose that $0<b<1$ and let  ${\bf u}(t,x;c)$ be a bounded classical  solution of \eqref{eq1} defined for every $x\in\mathbb{R}$ and $t\in\mathbb{R}$ such that $\min\{\inf_{x,t\in\mathbb{R}}u_2(t,x;c),\inf_{x,t\in\mathbb{R}}u_2(t,x)\}>0$. If there holds
\begin{equation}
\label{H6-eq}
 (1-\chi_1\mu_1M_1)_+(r-\chi_2\mu_2M_2)_+> abrM_1M_2,
 \end{equation}
 then ${\bf u}(t,x;c)\equiv {\bf e}_*$. In particular if ${\bf (H3)}$ holds then ${\bf u}(t,x;c)\equiv {\bf e}_*$.

\end{lemma}
\begin{proof}
First of all, note that {\bf (H3)} implies that
$$
1-\chi_1\mu_1M_1>aM_2\quad {\rm and}\quad r-\chi_2\mu_2M_2>brM_1,
$$
which implies \eqref{H6-eq}. It then suffices to prove ${\bf u}(t,x;c)\equiv {\bf e}_*$ provided that \eqref{H6-eq} holds.

Introducing $l_i$ and  $L_i$, $i=1,2$, as in \eqref{1-zzz-z11}, and noting that $\mu_1l_1+\mu_2l_2\leq\lambda v(x;{u})\leq \mu_1L_1+\mu_2L_2$ for every $x\in\mathbb{R}$, we can proceed as in the proof of Lemma \ref{lem-01} by using comparison principle for parabolic equations and the fact that $\min\{l_1,l_2,L_1,L_2\}>0$  to obtain the following inequalities :
\begin{equation}\label{1-k-1-1}
    L_1\leq M_1(1-\chi_1\mu_1l_1-al_2),
\end{equation}
\begin{equation}\label{1-k-2-1}
    l_1\geq M_1(1-\chi_1\mu_1L_1-aL_2),
\end{equation}
\begin{equation}\label{1-k-1-2}
  r  L_2\leq M_2(r-\chi_2\mu_2l_2-brl_1),
\end{equation}
\begin{equation}\label{1-k-2-2}
   r l_2\geq M_2(r-\chi_2\mu_2L_2-brL_1).
\end{equation}
Taking difference side by sides of inequalities \eqref{1-k-1-1} and \eqref{1-k-2-1} yields
\begin{equation}\label{1-k-3}
(1-\chi_1\mu_1M_1)(L_1-l_1)\leq aM_1(L_2-l_2).
\end{equation}
Similarly, it follows from inequalities \eqref{1-k-1-2} and \eqref{1-k-2-2} that
\begin{equation} \label{1-k-4}
(r-\chi_2\mu_2M_2)(L_2-l_2)\leq rbM_2(L_1-l_1).
\end{equation}
The last two inequalities imply that
$$
(1-\chi_1\mu_1M_1)(r-\chi_2\mu_2M_2)(L_1-l_1)(L_2-l_2)\leq abrM_1M_2(L_1-l_1)(L_2-l_2)
$$
since \eqref{H6-eq} implies that $1-\chi_1\mu_1M_1>0$ and $r-\chi_2\mu_2M_2>0$. Since \eqref{H6-eq} holds, we must have from the above inequality  that $(L_1-l_1)(L_2-l_2) =0$, which combined with \eqref{1-k-3} and \eqref{1-k-4} yield $l_1=L_1$ and $l_2=L_2$. Therefore, recalling the identities $(1-\chi_1\mu_1)M_1=1$ and $(r-\chi_2\mu_2)M_2=r$, it follows from \eqref{1-k-1-1}-\eqref{1-k-2-2} that
$$
\begin{cases}
L_1=1-aL_2\cr
L_2=1-bL_1.
\end{cases}
$$
It then follows that $l_1=L_1=\frac{1-a}{1-ab}$ and $l_2=L_2=\frac{1-b}{1-ab}$.  That is ${{\bf u}}(t,x;c)\equiv e_*$, which completes the proof of the lemma.
\end{proof}

Now, we prove Theorem \ref{stability-tm}.

\begin{proof}[Proof of Theorem \ref{stability-tm}] Let $c\in\mathbb{R}$ be given
and ${\bf u}(t,x;c)$ be a classical solution of \eqref{eq1} satisfying $\inf_{x\in\mathbb{R}}u_1(0,x;c)>0$. It  is easy to see that
$$
\limsup_{t\to\infty}\|u_i(t,\cdot;c)\|_{\infty}\le M_i \quad \forall i=1,2.
$$
Hence without loss of generality we may suppose that $\|u_i(t,\cdot;c)\|_{\infty}\le M_i$ for every $t\geq 0$ and $i=1,2$.

$(i)$  Suppose that hypothesis   ${\bf (H2)}$ holds. Observe that $u_1(t,x;c)$ satisfies
$$
    u_{1,t}\geq u_{1,xx}+(c-\chi_1v_{x})u_{1,x}+u_1(1-(\chi_1\mu_1M_1+aM_2)- (1-\chi_1\mu_1)u_1), \quad \forall\ t>0, \ x\in\mathbb{R}.
$$
Hence, since  ${\bf (H2)}$ holds,  we may employ the comparison principle for parabolic equations to conclude that
\begin{align}\label{am-1}
u_1(t,x;c)\ge \underline{m}_1:= \min\left\{\inf_{x\in\mathbb{R}}u_1(0,x;c),\frac{1-(\chi_1\mu_1M_1+aM_2)}{1-\chi_1\mu_1}  \right\}>0,\ \ \forall\ t\geq 0, \ x\in\mathbb{R}.
\end{align}\label{xdx-1}
Now, if we suppose by contradiction that the statement of Theorem \ref{stability-tm} $(i)$ is false, then there  exist a sequence $t_n\to\infty$ and $x_n$ such that
\begin{equation}\label{am-2}
    \inf_{n\ge 1}|{\bf u}(t_n,x_n;c)-{\bf e}_1|>0.
\end{equation}
By a priori estimates for parabolic equations, without loss of generality, we may suppose that there is some ${\bf u}^*\in C^{1,2}(\mathbb{R}\times\mathbb{R})$ such that ${\bf u}(t+t_n,x+x_n;c)\to{\bf u}^*(t,x)$ as $n\to\infty$ in $C^{1,2}_{\rm loc}(\mathbb{R}\times\mathbb{R})$. Note that ${\bf u}^*(t,x)$ is an entire solution of \eqref{eq1} and it follows from \eqref{am-1} that
$$
u_1(t,x)\geq \underline{m}_1,\quad \forall\ t,x\in\mathbb{R}.
$$
Thus, by Lemma \ref{lem-01}, we conclude that ${\bf u}^*(t,x)\equiv {\bf e}_1$. This contradicts with \eqref{am-2} since ${\bf u}^*(0,0)=\lim_{n\to\infty}{\bf u}(t_n,x_n;c)$. Hence $\lim_{t\to\infty}\|{\bf u}(t,\cdot;c)-{\bf e}_1\|_{\infty}=0$, which completes the proof of $(i)$.

\smallskip

$(ii)$  Suppose that hypothesis   ${\bf (H3)}$ holds. The proof follows similar arguments as in $(i)$. Indeed, note that in addition to \eqref{am-1}, it also holds that
$$
    u_{2,t}\geq du_{2,xx}+(c-\chi_2v_{x})u_{2,x}+u_2(r-(\chi_2\mu_2M_2+brM_1)- (r-\chi_2\mu_2)u_2), \quad \forall\ t>0, \ x\in\mathbb{R}.
$$
Hence, since  ${\bf (H3)}$ holds,  we may employ the comparison principle for parabolic equations to conclude that
\begin{align}\label{am-3}
u_2(t,x;c)\ge \underline{m}_2:= \min\left\{\inf_{x\in\mathbb{R}}u_2(0,x;c),\frac{r-(\chi_2\mu_2M_2+brM_1)}{r-\chi_2\mu_2}  \right\}>0,\ \ \forall\ t\geq 0, \ x\in\mathbb{R}.
\end{align}
Therefore, similar arguments used to prove $(i)$ together with Lemma \ref{lem-02} yield that $\lim_{t\to\infty}\|{\bf u}(t,\cdot;c)-{\bf e}_*\|_{\infty}=0$. This completes the proof of the theorem.
\end{proof}

\section{Super- and sub-solutions}\label{super-sol-sec}

In this section, we construct some  super- and sub-solutions  to some elliptic equations related to \eqref{eq1}. These super- and sub-solutions will be used in the proof of Theorem \ref{main-exist-thm} in next section.

 We first introduce some notations.  For fixed $0<\kappa<\kappa_{\max}=\min\{\sqrt{1-a},\sqrt{\lambda}\}$,  $D_1,D_2,\tilde{D}_2>0$  and $0<\varepsilon_1\ll 1$, we define
$$
\overline{u}^\kappa_1(x)=\min\{M_1,M_1D_2e^{-\kappa x}\}, \quad \underline{u}^\kappa_{1}(x)=M_1D_2\left(1-D_1e^{-\varepsilon_1 x} \right)_+e^{-\kappa x},
$$
$$
\overline{u}^\kappa_{2}(x)=\min\{M_2,1+M_2{\tilde{D}_{2}}e^{-\kappa x}\} \quad \text{and} \quad \underline{u}^\kappa_{2}(x)=\left(1-{\tilde{D}_{2}}e^{-\kappa x} \right)_+,
$$
where $(m)_+=\max\{m,0\}$ for every real number $m\in\mathbb{R}$, and $M_1$ and $M_2$ are as in \eqref{M1-M2-eq}, that is,  $M_1=\frac{1}{1-\chi_1\mu_1}$ and $M_2=\frac{r}{r-\chi_2\mu_2}$. We shall provide more information on how to choose the positive constants $D_1,D_2,{\tilde{D}_{2}}$ and $\varepsilon_1$ in Lemma \ref{lem3} below.
Define the convex set
\begin{equation}\label{set-E-def}
\mathcal{E}(\kappa):=\{{\bf u}\in C^b_{\rm unif}(\mathbb{R})\times C^b_{\rm unif}(\mathbb{R})\ |\  \ \underline{u}^\kappa_{i}(x)\leq u_i(x)\leq  \overline{u}^\kappa_{i}(x)\ \forall\ x\in\mathbb{R}, \ i=1,2 \}.
\end{equation}

Observe that for every ${\bf u}\in C^b_{\rm unif}(\mathbb{R})\times  C^b_{\rm unif}(\mathbb{R})$, the scalar valued function
\begin{equation}\label{eq2}
    v(x;{\bf u})=\frac{1}{2\sqrt{\lambda}}\int_{\mathbb{R}}e^{-\sqrt{\lambda}|x-y|}(\mu_1u_1(y)+\mu_2u_2(y))dy,\quad x\in\mathbb{R},
\end{equation}
is twice continuously differentiable   and solves the elliptic equation
$$
0=v_{xx}-\lambda v_x+\mu_1u_1+\mu_2u_2,\quad x\in\mathbb{R}.
$$

Next, we present some lemmas.

\medskip

\begin{lemma}\label{lem1}
Let ${\bf u}\in\mathcal{E}(\kappa)$ and $v(x;{\bf u})$ be given by \eqref{eq2}. Then for any $x\in\mathbb{R}$,
\begin{align}\label{xx-1}
    v(x;{\bf u})\leq \mu_1\min\{\frac{M_1}{\lambda},\frac{M_1D_2B_{\lambda,\kappa}}{2\sqrt{\lambda}}e^{-\kappa x}\}+\mu_2\min\{\frac{M_2}{\lambda},\frac{1}{\lambda}+\frac{M_2{\tilde{D}_{2}}B_{\lambda,\kappa}}{2\sqrt{\lambda}}e^{-\kappa x}\}
\end{align}
and
\begin{align}\label{xx-2}
v(x;{\bf u})\geq \frac{\mu_1M_1D_2 e^{-\kappa x}}{2\sqrt{\lambda}}\left(B_{\lambda,\kappa}-D_1B_{\lambda,\kappa+\varepsilon_1}e^{-\varepsilon_1 x}\right)_+ +\mu_2\left(\frac{1}{\lambda}-\frac{{\tilde{D}_{2}}B_{\lambda,\kappa}}{2\sqrt{\lambda}}e^{-\kappa x}\right)_+,
\end{align}
 where $B_{\lambda,\kappa}$ is as in \eqref{B-lambda-kappa-eq}.
\end{lemma}

\begin{lemma}\label{lem2}
Let ${\bf u}\in\mathcal{E}(\kappa)$ and $v(x;{\bf u})$ be given by \eqref{eq2}. Then
\begin{align}\label{xx-3}
   | \frac{d}{dx}v(x;{\bf u})|\leq  \mu_1\min\{\frac{M_1D_2B_{\lambda,\kappa}}{2}e^{-\kappa x},\frac{M_1}{\sqrt{\lambda}}\} +\mu_2\min\Big\{\frac{{\tilde{D}_{2}}M_2B_{\lambda,\kappa}}{2}e^{-\kappa x},\frac{M_2}{\sqrt{\lambda}}\Big\},\quad \forall\ x\in\mathbb{R},
\end{align}
where $B_{\lambda,\kappa}$ be as in \eqref{B-lambda-kappa-eq}.
\end{lemma}
We delay the proofs of Lemmas \ref{lem1} and \ref{lem2} to the Appendix. For every ${\bf u}\in\mathcal{E}(\kappa)$,  we associate the differential operator $\mathcal{F}^{\bf u}({\bf U})=(\mathcal{F}_{1}^{\bf u}( U_1),\mathcal{F}_{2}^{\bf u}( U_2) )$ defined by
\begin{equation}
    \mathcal{F}_{1}^{\bf u}(U_1)=U_{1,xx}+(c_\kappa-\chi_1v_x(\cdot;{\bf u}))U_{1,x}+(1-(a-\chi_1\mu_2)u_2-\lambda\chi_1v(\cdot;{\bf u})-\frac{(U_1)_+}{M_1})U_1+R(u_1-U_1),
\end{equation}
\begin{equation}
    \mathcal{F}_{2}^{\bf u}(U_2)=dU_{2,xx}+(c_\kappa-\chi_2v_x(\cdot,{\bf u}))U_{2,x}+(r-(br-\chi_2\mu_1)u_{1}-\lambda\chi_2v(\cdot;{\bf u})-\frac{r(U_2)_+}{M_2})U_2+R(u_2-U_2)
\end{equation}
with $(U_i)_+=\max\{U_i,0\}$, ${\bf U}=(U_1,U_2)$, where  $c_\kappa=\frac{\kappa^2+1-a}{\kappa}$ and $R\gg 1$ is a positive constant satisfying
\begin{equation}\label{eq:R-choice}
R>\max\{1-(a-\chi_1\mu_2)M_2-\chi_1(\mu_1M_1+\mu_2M_2),r-(br-\chi_2\mu_1)M_1-\chi_2(\mu_1M_1+\mu_2M_2)\}.
\end{equation}

Consider the elliptic system
\begin{equation}\label{comp-syst}
0=\mathcal{F}^{\bf u}({\bf U}),\quad x\in\mathbb{R}.
\end{equation}
  In the rest of this section, we assume that {\bf (H4)} holds.
 {Let $\kappa^*_{\chi_1,\chi_2}$ be defined as in \eqref{kappa-star-eq}. By the definition of $\kappa_{\chi_1,\chi_1}^*$, for every $0<\kappa<\kappa^*_{\chi_1,\chi_2}$,  it holds that
\begin{equation}\label{new-eq-2-1}
    \frac{\kappa\chi_1\mu_1M_1 B_{\lambda,\kappa}}{2}<1
\end{equation}  and
\begin{equation}\label{xxx-2}
   1-a\ge F(\kappa,\chi_1,\chi_2),
\end{equation}
where $F(\kappa,\chi_1,\chi_2)$ is as in \eqref{F-2-eq}.
} Our aim is to prove Theorem  \ref{main-exist-thm} in next section by showing that there is ${\bf u}^*\in\mathcal{E}(\kappa)$ such that $\mathcal{F}^{{\bf u}^*}({\bf u}^*)=0$
for $0<\kappa<\kappa^*:=\kappa^*_{\chi_1,\chi_2}$.

\begin{lemma}\label{lem3} Let {$0<\kappa<\kappa^*$},  $0<\varepsilon_1<\min\{\kappa,c_\kappa-2\kappa\},$ $D_2=\frac{1}{M_1}$,  and  $\tilde{D}_2=D_2/ (a+\chi_1f(\kappa,\chi_1,\chi_2))$, where $f(\kappa,\chi_1,\chi_2)$ is the positive solution of the algebraic equation  \eqref{new-f-de-1}.
Then there is $D_1\gg 1$ such that for every ${\bf u}\in \mathcal{E}(\kappa)$  and $i\in\{1,2\}$ the following hold.
\begin{itemize}
    \item [(i)] $\mathcal{F}_i^{\bf u}(M_i)\leq 0$ and $\mathcal{F}_i^{\bf u}(0)\ge 0$ for every $x\in\mathbb{R}$. 
    \item[(ii)] $\mathcal{F}_i^{\bf u}(\overline{u}^\kappa_i)\le 0$ on the open interval $\{\overline{u}^\kappa_i<M_i\}$. 
    \item[(iii)]  $\mathcal{F}_i^{\bf u}(\underline{u}^\kappa_i)\ge 0$ on the open interval $\{\underline{u}^\kappa_i>0\} $. 
\end{itemize}
\end{lemma}

\begin{proof} Let ${\bf u}\in\mathcal{E}(\kappa)$.

$(i)$ Note that $\mathcal{F}_i^{\bf u}(0)=Ru_i(x)\ge 0$ for every $x\in\mathbb{R}$ and
$$
\mathcal{F}_1^{\bf u}(M_1)=-M_1((a-\chi_1\mu_2)u_2(x)+\chi_1\lambda v(x;{\bf u}))-R(M_1-u_1(x))<0, \quad \forall\ x\in\mathbb{R}
$$
since $v(\cdot;{\bf u}) > 0$, $\|u_1\|_{\infty}\leq M_1$ and $u_2(\cdot)\ge 0$. Similarly
$$
\mathcal{F}_2^{\bf u}(M_2)=-M_2((br-\chi_2\mu_1)u_1(x)+\chi_2\lambda v(x;{\bf u}))-R(M_2-u_2(x))<0, \quad \forall\ x\in\mathbb{R}
$$
since $v(\cdot;{\bf u}) \geq 0$, $\|u_2\|_{\infty}\leq M_2$ and $u_1(\cdot)\ge 0$.

\smallskip

 $(ii)$ Recalling inequalities \eqref{xx-2} and \eqref{xx-3},  with $\overline{u}^\kappa_1=M_1D_2e^{-\kappa x} $, we get
 
\begin{align*}
    \mathcal{F}_1^{\bf u}(\overline{u}^\kappa_1)=&\left(\kappa^2-\kappa c_\kappa+1+\kappa\chi_1v_x(;{\bf u})-D_2e^{-\kappa x}-(a-\chi_1\mu_2)u_2-\chi_1\lambda v(\cdot;{\bf u})) \right)\overline{u}^\kappa_1 -R(\overline{u}^\kappa_1 -u_1)\cr
    = &\left(a+\kappa\chi_1v_x(;{\bf u})-D_2e^{-\kappa x}-(a-\chi_1\mu_2)u_2-\chi_1\lambda v(\cdot;{\bf u})) \right)\overline{u}^\kappa_1-R(\overline{u}^\kappa_1 -u_1)\cr
    \leq & \left(a-\chi_1\mu_2+\kappa\chi_1v_x(;{\bf u})-\tilde{D}_{2}\Big(\frac{D_2}{\tilde{D}_{2}}-\frac{\chi_1\mu_2\sqrt{\lambda} B_{\lambda,\kappa}}{2}\Big)e^{-\kappa x}-(a-\chi_1\mu_2)\underline{u}^\kappa_{2} \right)\overline{u}^\kappa_1\cr
    \leq & \left(a-\chi_1\mu_2+\kappa\chi_1v_x(;{\bf u})-\tilde{D}_2\Big(\frac{D_2}{\tilde{D}_{2}}-\frac{{D}_{2}\chi_1\mu_2\sqrt{\lambda} B_{\lambda,\kappa}}{2}\Big)e^{-\kappa x}-(a-\chi_1\mu_2)(1-\tilde{D}_{2}e^{-\kappa x}) \right)\overline{u}^\kappa_1\cr
    =& \left(\kappa\chi_1v_x(;{\bf u})-\tilde{D}_2\Big(\frac{D_2}{\tilde{D}_{2}}-\frac{\chi_1\mu_2\sqrt{\lambda} B_{\lambda,\kappa}}{2} -(a-\chi_1\mu_2)\Big)e^{-\kappa x} \right)\overline{u}^\kappa_1\cr
    \leq & \tilde{D}_2\left(\frac{\chi_1\kappa B_{\lambda,\kappa}}{2}(\frac{\mu_1M_1D_2}{\tilde{D}_2}+\mu_2M_2)-\Big(\frac{D_2}{\tilde{D}_2}-\frac{\chi_1\mu_2\sqrt{\lambda} B_{\lambda,\kappa}}{2}-(a-\chi_1\mu_2)\Big) \right)\overline{u}^\kappa_1e^{-\kappa x}\cr
    =&   \chi_1
    \tilde{D}_2\Big(\frac{\kappa B_{\lambda,\kappa}}{2}\big( \mu_1M_1(a+\chi_1f)+\mu_2M_2 \big) + \frac{\mu_2\kappa^2B_{\lambda,\kappa}}{2\sqrt{\lambda}}-f \Big)\overline{u}^\kappa_1 e^{-\kappa x} =0.
\end{align*}
Note that we have used \eqref{new-f-de-1}.
Similarly,  for $\overline{u}^\kappa_2(x)=1+M_2D_2e^{-\kappa x}\le M_2$ and using Lemmas \ref{lem1} and \ref{lem2}, we obtain

\begin{align*}
\begin{split}
    \mathcal{F}_2^{\bf u}(\overline{u}^\kappa_2)
    =& M_2(d\kappa^2-\kappa(c_\kappa-\chi_2v_x(\cdot;{\bf u})))\frac{\tilde{D_2}}{e^{\kappa x}}+(r-\frac{r\overline{u}^\kappa_2}{M_2}-(br-\chi_2\mu_1)u_1-\chi_2\lambda v(\cdot;{\bf u}))\overline{u}^\kappa_2+R(u_2-\overline{u}^\kappa_2)\cr
    \le & M_2(d\kappa^2-\kappa(c_\kappa-\chi_2v_x(\cdot;{\bf u})))\tilde{D_2}e^{-\kappa x} +(r-\frac{r}{M_2}\overline{u}^\kappa_2-\chi_2\lambda v(\cdot;{\bf u}))\cr
    = & M_2(d\kappa^2-\kappa(c_\kappa-\chi_2v_x(\cdot;{\bf u})))\tilde{D_2}e^{-\kappa x}+\Big(r-\frac{r}{M_2}(1+M_2\tilde{D_2}e^{-\kappa x})-\chi_2\lambda v(\cdot;{\bf u}) \Big)\overline{u}^\kappa_2\cr
    \le &M_2(d\kappa^2-\kappa(c_\kappa-\chi_2v_x(\cdot;{\bf u})))\tilde{D_2}e^{-\kappa x}+ \Big(r-\frac{r}{M_2}(1+M_2\tilde{D_2}e^{-\kappa x}) -\chi_2\mu_2\Big(1-\frac{\tilde{D_2}\sqrt{\lambda}B_{\lambda,\kappa}}{2e^{\kappa x}} \Big)  \Big)\overline{u}^\kappa_2\cr
    = & \tilde{D_2}\left((d\kappa^2-\kappa(c_\kappa-\chi_2v_x(\cdot;{\bf u})))M_2-\Big(r-\frac{\chi_2\mu_2\sqrt{\lambda}B_{\lambda,\kappa}}{2}  \Big)\overline{u}^\kappa_2\right)e^{-\kappa x}\cr
    \leq& \tilde{D_2}\left((d\kappa^2-\kappa c_\kappa+\frac{\kappa\chi_2}{\sqrt{\lambda}}(\mu_1M_1+\mu_2M_2))M_2-\Big(r-\frac{\chi_2\mu_2\sqrt{\lambda}B_{\lambda,\kappa}}{2}  \Big)\overline{u}^\kappa_2\right)e^{-\kappa x}\cr
    \leq & \frac{M_2\tilde{D_2}}{e^{\kappa x}}\left((d-1)\kappa^2-(1-a)+\frac{\kappa\chi_2}{\sqrt{\lambda}}(\mu_1M_1+\mu_2M_2)+\frac{\chi_2\mu_2\sqrt{\lambda}B_{\lambda,\kappa}}{2}\frac{\overline{u}^\kappa_2}{M_2}-\frac{r}{M_2}\right)\cr
    \le & \frac{M_2\tilde{D_2}}{e^{\kappa x}}\left((d-1)\kappa^2-(1-a)+\frac{\kappa\chi_2}{\sqrt{\lambda}}(\mu_1M_1+\mu_2M_2)+\frac{\chi_2\mu_2\sqrt{\lambda}B_{\lambda,\kappa}}{2}  -\frac{r}{M_2}\right)\le 0
    \end{split}
\end{align*}
whenever \eqref{xxx-2} holds. This completes the proof of $(ii)$.

\smallskip

$(iii)$ For $x\in\{\underline{u}^\kappa_{2}>0\}$ and recalling Lemmas \ref{lem1} and \ref{lem2}, we obtain
\begin{small}
\begin{align*}
&  \mathcal{F}_2^{\bf u}(\underline{u}^\kappa_{2})+R(\underline{u}^\kappa_{2}-u_2)\cr
  =&\left(r-\frac{r}{M_2}(1-\tilde{D_2}e^{-\kappa x})-(br-\chi_2\mu_1)u_1-\chi_2\lambda v(\cdot;{\bf u}) \right)\underline{u}^\kappa_{2} + (\kappa c_\kappa-d\kappa^2-\kappa\chi_2v_x(\cdot;{\bf u}))\tilde{D_2}e^{-\kappa x}  \cr
  =& \left(r-\frac{r}{M_2}(1-\tilde{D_2}e^{-\kappa x})-(br-\chi_2\mu_1)u_1-\chi_2\lambda v(\cdot;{\bf u}) \right)\underline{u}^\kappa_{2} + (1-a+(1-d)\kappa^2-\kappa\chi_2v_x(\cdot;{\bf u}))\tilde{D_2}e^{-\kappa x}  \cr
\geq& \left(r-\frac{r}{M_2}(1-\tilde{D_2}e^{-\kappa x})-(br-\chi_2\mu_1)u_1-\chi_2\lambda v(\cdot;{\bf u}) \right)\underline{u}^\kappa_{2} + \Big(1-a+(1-d)\kappa^2
  -\frac{\kappa\chi_2(\mu_1 M_1+\mu_2M_2)}{\sqrt{\lambda}}\Big)\frac{\tilde{D_2}}{e^{\kappa x } } \cr
  \ge & \left(r-\frac{r(1-\tilde{D_2}e^{-\kappa x})}{M_2}-\Big(\frac{M_1D_2}{\tilde{D}_2}(br-\chi_2\mu_1)+\frac{\chi_2\sqrt{\lambda} B_{\lambda,\kappa}( \mu_1M_1+\mu_2M_2)}{2}\Big)\frac{\tilde{D_2}}{e^{\kappa x}}-\chi_2\mu_2\right)\underline{u}^\kappa_{2}\cr
  & + \Big(1-a+(1-d)\kappa^2-\frac{\kappa\chi_2(\mu_1 M_1+\mu_2M_2)}{\sqrt{\lambda}}\Big)\tilde{D_2}e^{-\kappa x}  \cr
  =&-\tilde{D_2}\left(r\big(M_1(a+\chi_1f(\kappa,\chi_1,\chi_2))(b-\frac{\chi_2\mu_1}{r})-\frac{1}{M_2}\big)+\frac{\chi_2 \sqrt{\lambda}B_{\lambda,\kappa}( \mu_1M_1+\mu_2M_2)}{2}\right)\underline{u}^\kappa_{2}e^{-\kappa x}\cr
  & + \Big(1-a+(1-d)\kappa^2 -\frac{\kappa\chi_2(\mu_1 M_1+\mu_2M_2)}{\sqrt{\lambda}}
  \Big)\tilde{D_2}e^{-\kappa x}
 \cr
 \geq& \left(1-a-(d-1)\kappa^2-r\left(M_1\big(a+\chi_1f(\kappa,\chi_1,\chi_2)\big)\big(b-\frac{\chi_2\mu_1}{r}\big)-\frac{1}{M_2}\right)_{+}-\frac{\chi_2(\lambda B_{\lambda,\kappa}+2\kappa)( \mu_1M_1+\mu_2M_2)}{2\sqrt{\lambda}}\right)\frac{\tilde{D_2}}{e^{\kappa x}} \cr
 =& ((1-a)-F(\kappa,\chi_1,\chi_2))\tilde{D_2}e^{-\kappa x}
 \ge  0
\end{align*}
\end{small}
 since \eqref{xxx-2} holds.
On the other hand, for $x\in \{\underline{u}^\kappa_{1} >0\}$, it holds
\begin{align*}
    \mathcal{F}_1^{\bf u}(\underline{u}^\kappa_{1})= &D_2M_1\left(\kappa^2-D_1(\varepsilon_1+\kappa)^2e^{-\varepsilon_1 x} \right)e^{-\kappa x}+D_2M_1\left( D_1(\varepsilon_1+\kappa)e^{-\varepsilon_1 x}-\kappa\right)(c_\kappa-\chi_1v_{x}(\cdot;{\bf u}))e^{-\kappa x}\cr
    & +(1-\frac{1}{M_1}\underline{u}^\kappa_{1}-(a-\chi_1\mu_2)u_2 -\chi_1\lambda v(\cdot;{\bf u}) )\underline{u}^\kappa_{1} +R(u_1-\underline{u}^\kappa_{1})\cr
    =& (a-1)\underline{u}^\kappa_{1}+D_2M_1e^{-\kappa x}\left( D_1\varepsilon_1(c_\kappa-\varepsilon_1-2\kappa)e^{-\varepsilon_1 x}+\chi_1(\kappa-D_1(\varepsilon_1+\kappa)e^{-\varepsilon_1 x})v_x(\cdot;{\bf u})\right)\cr
    & +(1-\frac{1}{M_1}\underline{u}^\kappa_{1}-(a-\chi_1\mu_2)u_2 -\chi_1\lambda v(\cdot;{\bf u}) )\underline{u}^\kappa_{1} +R(u_1-\underline{u}^\kappa_{1})\cr
    \ge &  (a-1)\underline{u}^\kappa_{1}+D_2M_1e^{-\kappa x}\left( D_1\varepsilon_1(c_\kappa-\varepsilon_1-2\kappa)e^{-\varepsilon_1 x}+\chi_1(\kappa-D_1(\varepsilon_1+\kappa)e^{-\varepsilon_1 x})v_x(\cdot;{\bf u})\right)\cr
    +& \Big(1-\frac{1}{M_1}\underline{u}^\kappa_{1}-(a-\chi_1\mu_2)(1+M_2\tilde{D}_{2}e^{-\kappa x}) -\chi_1\mu_2- \frac{D_2\chi_1}{2}(\mu_2M_2\frac{\tilde{D}_2}{D_2}+\mu_1M_1)\sqrt{\lambda}B_{\lambda,\kappa}e^{-\kappa x} \Big)\underline{u}^\kappa_{1} \cr
   = & D_2M_1e^{-\kappa x}\left( D_1\varepsilon_1(c_\kappa-\varepsilon_1-2\kappa)e^{-\varepsilon_1 x}+\chi_1(\kappa-D_1(\varepsilon_1+\kappa)e^{-\varepsilon_1 x})v_x(\cdot;{\bf u})\right)\cr
    & -\Big(\frac{1}{M_1}\underline{u}^\kappa_{1}+ D_2\left(  \frac{(a-\chi_1\mu_2)M_2\tilde{D}_2}{D_2}+\frac{\chi_1}{2}(\frac{\mu_2M_2\tilde{D}_2}{D_2}+\mu_1M_1)\sqrt{\lambda}B_{\lambda,\kappa} \right)e^{-\kappa x} \Big)\underline{u}^\kappa_{1} \cr
    \geq& D_2M_1e^{-\kappa x}\left( D_1\varepsilon_1(c_\kappa-\varepsilon_1-2\kappa)e^{-\varepsilon_1 x}+\chi_1(\kappa-D_1(\varepsilon_1+\kappa)e^{-\varepsilon_1 x})v_x(\cdot;{\bf u})\right)\cr
    & -D_2\Big(1+   \frac{(a-\chi_1\mu_2)M_2\tilde{D}_2}{D_2}+\frac{\chi_1}{2}(\frac{\mu_2M_2\tilde{D}_2}{D_2}+\mu_1M_1)\sqrt{\lambda}B_{\lambda,\kappa} \Big)M_1D_2e^{-2\kappa x}
\end{align*}
where we have used $\underline{u}^\kappa_{1}(x)\leq D_2M_1e^{-\kappa x}$ and Lemma \ref{lem1}. Using Lemma \ref{lem2} and the fact that $D_1e^{-\varepsilon_1 x}\leq 1$ and $e^{-\varepsilon_1 x}>e^{-\kappa x}$ for $x\in\{\underline{u}^\kappa_{1}>0\}$ and $D_1>1$, the last inequality is improved to
\begin{align}\label{xxx-1}
    \mathcal{F}_1^{\bf u}(\underline{u}^\kappa_{1})\geq& D_2M_1e^{-2\kappa x}\left( D_1\varepsilon_1(c_\kappa-\varepsilon_1-2\kappa)-\frac{D_2\chi_1(2\kappa+\varepsilon_1)}{2}(\mu_1M_1+\frac{\mu_2M_2\tilde{D}_2}{D_2})B_{\lambda,\kappa}\right)\cr
    & -D_2\Big(1+   \frac{(a-\chi_1\mu_2)M_2\tilde{D}_2}{D_2}+\frac{\chi_1}{2}(\frac{\mu_2M_2\tilde{D}_2}{D_2}+\mu_1M_1)\sqrt{\lambda}B_{\lambda,\kappa} \Big)D_2M_1e^{-2\kappa x}.
\end{align}
Hence with the choice of $\frac{D_1}{D_2}=D_1M_1\gg 1$ satisfying
\begin{equation}
    \frac{D_1}{D_2}\varepsilon_1(c_\kappa-\varepsilon_1-2\kappa)\ge 1+   \frac{(a-\chi_1\mu_2)M_2\tilde{D}_2}{D_2}+\frac{\chi_1(\sqrt{\lambda}+2\kappa+\varepsilon_1)}{2}(\frac{\mu_2M_2\tilde{D}_2}{D_2}+\mu_1M_1)B_{\lambda,\kappa}
\end{equation} 
we conclude from \eqref{xxx-1} that
$$
\mathcal{F}_1^{\bf u}(\underline{u}^\kappa_{1})(x)\ge 0 \quad \forall\ x\in\{\underline{u}^\kappa_{1}>0\},
$$
which completes the proof of $(iii)$.
\end{proof}

\section{Existence of traveling wave solutions for $c>c^*$}\label{sec-for-t1}

In this section, we investigate the existence of traveling wave solutions of \eqref{eq:1} and prove  Theorem \ref{main-exist-thm}.

 {Throughout this section, we assume that {\bf (H4)} holds and that  the constants $\varepsilon_1,D_2,\tilde{D}_2$ and $D_1$ are  fixed and satisfy the hypotheses of Lemma \ref{lem3}.} Recall that $c^*=\frac{(k^*)^2+1-a}{k^*}$, where $\kappa^*={\kappa^*_{\chi_1,\chi_2}}$
is as in \eqref{kappa-star-eq}.
  Our idea to prove Theorem \ref{main-exist-thm}  is to prove that, for any $0<\kappa<\kappa^*$,  there is ${\bf u}^*\in\mathcal{E}(\kappa)$ such that $\mathcal{F}^{{\bf u}^*}({\bf u}^*)=0$.

\smallskip

To this end,
for every $y>0 $  {and $0<\kappa<\kappa^*$}, and ${\bf u}\in \mathcal{E}(\kappa)$, consider the following elliptic system
\begin{equation}\label{y-dependent-equation}
    \begin{cases}
    0=\mathcal{F}_i^{{\bf u}}({ U}_i^{{\bf u},y}),& |x|<y,\ i=1,2\cr
    {{ U}_i^{{\bf u},y}}( x)=\overline{u}^\kappa_i(x), & x=\pm y,\ i=1,2.
    \end{cases}
\end{equation}


\begin{lemma}\label{lem4}
{For every $y>0$ and  $0<\kappa<\kappa^*$}, and ${\bf u}\in \mathcal{E}(\kappa)$ there exists a unique ${{\bf U}^{{\bf u},y}}$ satisfying \eqref{y-dependent-equation}.
\end{lemma}
\begin{proof} Let ${\bf u}\in\mathcal{E}(\kappa)$ and $y>0$ be given.
We first show the uniqueness. Observe that system \eqref{y-dependent-equation} is decoupled, hence the theory of elliptic scalar equations applies for each equation. Since the equations of \eqref{y-dependent-equation} are of the same type for both $i=1$ and $i=2$, we shall only provide the arguments for the proof of $U_1^{{\bf u},y}$.

Observe from the choice of the positive constant $R$, see \eqref{eq:R-choice}, and Lemma \ref{lem1}  that
$$
1-(a-\chi_1\mu_1)u_2(x)-\lambda\chi_1v(x;{\bf u})< R , \quad\forall\ x\in\mathbb{R}.
$$
Whence for each $x\in\mathbb{R}$ fixed, the function
\begin{equation}\label{reaction-term-1}
\mathbb{R}\ni U_1\mapsto \mathcal{A}_1^{{\bf u}}(x,U_1):=(1-(a-\chi_1\mu_2)u_{2}(x)-\lambda\chi_1v(x,{\bf u})-\frac{(U_1)_+}{M_1})U_1+R(u_1(x)-U_1)
\end{equation}
is monotone decreasing. Thus by \cite[Theorem 10.2, page 264 ]{Trudinger} we deduce that a solution $U^{{\bf u},y}(x)$, if exists, is unique.

 Now, we show the existence of solution to \eqref{y-dependent-equation}. Again, we note from the choice of $R$ (see \eqref{eq:R-choice}) and Lemma \ref{lem1} that
 $$
  \mathcal{A}_1^{{\bf u}}(x,M_1)\le 0 \quad \text{and}\quad \mathcal{A}_1^{{\bf u}}(x,-M_1) \geq 0
 $$
 for every $x\in\mathbb{R}$. We also note, by Lemma \ref{lem2},
 $$
 |c_{\kappa}-\chi_1v_x(x;{\bf u})|\leq c_{\kappa}+\frac{\chi_1}{\sqrt{\lambda}}\left(\mu_1M_1+\mu_2M_2 \right),\quad \forall \ x\in\mathbb{R}.
 $$
 Since $\|\overline{u}^\kappa_1\|_{\infty}\leq M_1$, it follows from \cite[Theorem 5.1, Corollary 5.2, page 433]{Hartman} that there is at least one classical solution to the elliptic equation
 $$
 \begin{cases}
 0=\mathcal{F}_1^{{\bf u},y}(U_1) & |x|<y\cr
 U_1(x)=\overline{u}^\kappa_1(x) & x=\pm y.
 \end{cases}
 $$

\end{proof}

For reference later, we  introduce the function
\begin{equation}\label{reaction-term-2}
    \mathcal{A}_2^{{\bf u}}(x,U_2)=(r-(br-\chi_2\mu_1)u_1-\lambda\chi_2v(x,{\bf u})-\frac{r}{M_2}(U_2)_+)U_2 +R(u_2(x)-U_2).
\end{equation}
Then \eqref{y-dependent-equation} can be rewritten as
\begin{equation}\label{rewritten-eq}
\begin{cases}
    0=d_iU^{{\bf u},y}_{i,xx}+(c_\kappa-\chi_iv_x(\cdot;{\bf u}))U_{i,xx}^{{\bf u},y} +\mathcal{A}_i^{{\bf u}}(U_i^{{\bf u},y}) & |x|<y\cr
    U^{{\bf u},y}_i(x)=\overline{u}^\kappa_i(x) & x=\pm y
    \end{cases}
\end{equation}
with $d_1=1$ and $d_2=d$.  {For every $y>0$,    $0<\kappa<\kappa^*$}, and ${\bf u}\in \mathcal{E}(\kappa)$  we define ${{\bf U}^{{\bf u},y}}(x)=(\overline{u}^\kappa_{1}(x),\overline{u}^\kappa_{2}(x))$  for every $|x|>y$. With this extension, we have the following lemma. For convenience we let $y_0>1$ such that $\overline{u}^\kappa_i(-y)=M_i$ for each $i\in\{1,2\}$ and $y\ge y_0$.

\begin{lemma}\label{lem5}
Let {$0<\kappa<\kappa^*$}, ${\bf u}\in\mathcal{E}(\kappa)$,  and $y\ge y_0$ be given. The following hold for every $i\in\{1,2\}$.
\begin{itemize}
    \item[(i)] For every $x\in\mathbb{R}$, $\underline{u}^\kappa_{i}(x)\leq { U}_i^{{\bf u},y}(x)\leq \overline{u}^\kappa_{i}(x) $
    \item[(ii)] If ${\bf (H2)}$ holds and $U_1^{{\bf u},y}(x)\equiv u_1(x)$, then there exist $0<m_1^*\ll 1$  and $x_1>0$ such that $U_1^{{\bf u},y}(x)\geq m_1^*$ for every $-y\leq x\leq x_1$ whenever $y>x_1$.
    \item[(iii)] If  {\bf (H3)} holds, and  $U_2^{{\bf u},y}(x)\equiv u_2(x)$, then there exist $0< m_2\ll 1$ and $x_2>0$ such that $U_2^{{\bf u},y}(x)\geq m_2$ for every $-y\leq x\leq x_2$ whenever  $y>x_2$.
\end{itemize}
\end{lemma}
\begin{proof}
$(i)$ Let $i\in\{1,2\}$. Since for every $x\in\mathbb{R}$ fixed, the function $U_i\mapsto \mathcal{A}_i^{\bf u}(U_i)$ is monotone decreasing, then it follows from Lemma \ref{lem3}  and  the comparison principle for elliptic equations (see \cite[Theorem 10.1, page 263]{Trudinger}) that
$$
\underline{u}^\kappa_i(x)\leq U^{{\bf u},y}_i(x)\le \overline{u}^\kappa_i(x),\quad \forall\ -y<x<y,
$$
which together with the fact that $U^{{\bf u},y}(x)=\overline{u}^\kappa_1(x)$ for every $|x|\geq y$, completes the proof of  $(i)$.

\smallskip

$(ii)$ Suppose that ${\bf (H2)}$ holds and that $U_1^{{\bf u},y}(x)= u_1(x)$ for every $|x|\le y$. Hence $u_1(x)$ satisfies
\begin{equation}\label{hhh-1}
    \begin{cases}
    0=u_{1,xx}+(c_\kappa-\chi_1v_x(\cdot;{\bf u}))u_{1,x}+(1-(1-\chi_1\mu_1)u_1-(a-\chi_1\mu_2)u_2-\chi_1\lambda v(\cdot;{\bf u}))u_1  & |x|<y\cr
    u_1(x)=\overline{u}^{\kappa}_1(x) & x=\pm y.
    \end{cases}
\end{equation}
Since $u_1(x)\ge \underline{u}_1^{\kappa}(x)\ge0$ and $u_1(\pm y)>0$, then the Harnack's inequality for elliptic equations implies that $u_1(x)>0$ for every $|x|\le y$. Observe that with $x_1:=\frac{1}{\varepsilon_1}\ln\Big(\frac{D_1(\varepsilon_1+\kappa)}{\kappa}\Big)$, it holds that
$$
0<m_1:=\max_{x\in\mathbb{R}}\underline{u}^\kappa_1(x)=\underline{u}^\kappa_{1}(x_1)=\frac{D_2M_1\varepsilon_1\kappa^{\frac{\kappa}{\varepsilon_1}}}{(\kappa+\varepsilon_1)^{\frac{\kappa}{\varepsilon}+1}}D_1^{-\frac{\kappa}{\varepsilon_1}}=\frac{\varepsilon_1\kappa^{\frac{\kappa}{\varepsilon_1}}}{(\kappa+\varepsilon_1)^{\frac{\kappa}{\varepsilon}+1}}D_1^{-\frac{\kappa}{\varepsilon_1}}\le u_1(x_1)
$$
and that
$$
1>\eta_1:=1-(a-\chi_1\mu_2)M_2-\chi_1(\mu_1M_1-\mu_2M_2)=1-aM_2-\chi_1\mu_1M_1>0.
$$
Now we take $m_1^*:=\min\{\eta_1M_1,m_1\}$. We claim that for every $y>\max\{y_0,x_1\}$ it holds that $m_1^*\leq u_1(x)$ for every $x\in[-y,x_1]$. Indeed, let $y>\max\{x_1,y_0\}$ and suppose that $u_1(x)$ attains its minimum at some point, say $\tilde{x}_1\in[-y,x_1]$. If $\tilde{x}_1$ is a boundary point, then the claim easily follows. On the other hand, if $\tilde{x}_1$ is an interior point then $u_{1,x}(\tilde{x}_1)=0$ and $u_{1,xx}(\tilde{x}_1)\ge0$. This along with \eqref{hhh-1} and the fact that $0\leq u_i\le M_i$ for each $i=1,2$ imply that
\begin{align*}
0\geq (1-(1-\chi_1\mu_1)u_1(\tilde{x}_1)-(a-\chi_1\mu_2)u_2(\tilde{x}_1)-\chi_1\lambda v((\tilde{x}_1);{\bf u}))u_1(\tilde{x}_1)
\geq  (\eta_1-(1-\chi_1\mu_1)u_1(\tilde{x}_1))u_1(\tilde{x}_1).
\end{align*}
This clearly implies that $u_1(\tilde{x}_1)\geq \frac{\eta_1}{1-\chi_1\mu_1}=\eta_1M_1\ge m^*_1$ since we have shown in the above that $u_1(\tilde{x}_1)=\min_{|x|\le y }u_1(x)>0$. So, the claim holds and the result follows.


\smallskip

$(iii)$ Suppose that ${\bf (H3)}$ holds and $U_2^{\bf u,y}(x)=u_2(x)$ for $|x|\le y$. The proof follows similar arguments as in $(ii)$ by observing that
$$
\eta_2:=r-(br-\chi_2\mu_1)M_1-\chi_2(\mu_1M_1+\mu_2M_2)=r-rbM_1-\chi_2\mu_2M_2>0
$$
and  with $x_2=\frac{\ln(\tilde{D}_{2}+\varepsilon_2)}{\kappa}$,
$$
u_2(x_2)>m_2=\underline{u}^\kappa_2(x_2)=1-\frac{\tilde{D}_2}{\tilde{D}_2+\varepsilon_2}=\frac{\varepsilon_2}{\tilde{D}_2+\varepsilon_2}\to 0 \quad \text{as}\ \varepsilon_2\to0^+.
$$
Hence we can take $0<\varepsilon_2\ll 1$ such that  we take $m_2^*=m_2=\min\{m_2,\eta_2M_2,M_2\}$ and then proceed as in the proof of $(ii)$ to show that $u_2(x)\geq m_2^*$ for every $x\in[-y,x_2]$ with $y>\max\{x_2,y_0\}$.
\end{proof}

\begin{theorem}\label{tm-1}
Let {$0<\kappa<\kappa^*$}  and $y\ge y_0$ be given. Then ${{\bf U}^{{\bf u},y}}\in  \mathcal{E}(\kappa)$ for every  ${\bf u}\in \mathcal{E}(\kappa)$. Moreover, the mapping $  \mathcal{E}(\kappa)\ni {\bf u} \mapsto {{\bf U}^{{\bf u},y}}\in  \mathcal{E}(\kappa) $ is continuous and compact with respect to the compact open topology. Therefore, by  Schauder's fixed point theorem, it has a fixed point, say ${\bf u}^{*,y}$.
\end{theorem}
\begin{proof}
Let $y>y_0$ be fixed. It follows from Lemma \ref{lem5} (i) that $ {\bf U}^{{\bf u},y}\in\mathcal{E}(\kappa)$ for every ${\bf u}\in\mathcal{E}(\kappa)$. Since $\|U^{{\bf u},y}_i\|_{\infty}\leq M_i$ for every $i=1,2$ and $u\in\mathcal{E}(\kappa)$, by a priori estimates for elliptic equations and the uniqueness of solution to \eqref{y-dependent-equation} guaranteed by Lemma \ref{lem4} and the Arzela-Ascot's theorem, it follows that  the mapping $  \mathcal{E}(\kappa)\ni {\bf u} \mapsto {\bf U}^{{\bf u},y}\in  \mathcal{E}(\kappa) $ is continuous and compact with respect to the compact open topology. Therefore, by  Schauder's fixed point theorem, it has a fixed point, say ${\bf u}^{*,y}$.
\end{proof}

{Fix $0<\kappa<\kappa^*$}.
For every $y>y_0$,  let ${\bf u}^{*,y}$ be a fixed point of the mapping $  \mathcal{E}(\kappa)\ni {\bf u} \mapsto {\bf U}^{{\bf u},y}\in  \mathcal{E}(\kappa) $. Since $\|u^{*,y}_{i}\|_{\infty}\leq M_1+M_2$, $i=1,2$, and by Lemmas \ref{lem1} and \ref{lem2} it holds that
$$
\|v(\cdot;{\bf u}^{*,y})\|_{C^{2,b}_{\rm unif}(\mathbb{R})}\leq \lambda\|v(\cdot;{\bf u}^{*,y})\|+\mu_1\|u^{*,y}_{1}\|_{\infty}+\mu_2\|u^{*,y}_{2}\|_{\infty}\leq 2(\mu_1M_1+\mu_2M_2)
$$
for every $y>y_0$. By a priori estimates for elliptic equations and  Arzela-Ascoli's theorem, there is a sequence $\{y_n\}_{n\geq 1}$  with $y_n\to\infty$ as $n\to\infty$ and a function ${\bf u}^*\in C^{2,b}_{\rm unif}(\mathbb{R})\times  C^{2,b}_{\rm unif}(\mathbb{R})$ such that ${\bf u}^{*,y_n}\to {\bf u}^*$ locally uniformly in $C^2(\mathbb{R})\times C^2(\mathbb{R})$. Moreover, the function ${\bf u}^*$  satisfies the elliptic system
\begin{equation}\label{zzz-z2}
    \begin{cases}
    0={u}_{1,xx}^*+(c_\kappa-\chi_1v^*_x(\cdot)) u^*_{1,x}+(1-(1-\chi_1\mu_1)u^*_1-(a-\chi_1\mu_2)u^*_2-\chi_1\lambda v^*)u^*_1, & x\in\mathbb{R}\cr
    0=d{u}_{1,xx}^*+(c_\kappa-\chi_2v^*_x(\cdot)) u^*_{1,x}+(r-(r-\chi_2\mu_2)u^*_2-(br-\chi_2\mu_1)u^*_1-\chi_2\lambda v)u^*_1, & x\in\mathbb{R}\cr
    0=v^*_{xx}-\lambda v^*+\mu_1 u_1^*+\mu_2 u^*_2, & x\in\mathbb{R}
    \end{cases}
\end{equation}
where $v^*(\cdot)=v(\cdot;{\bf u}^*)$. Since ${\bf u}^*\in\mathcal{E}(\kappa)$, then
\begin{equation}\label{zzz-z3-1-1-1}
\lim_{x\to\infty}\frac{u^*_1(x)}{e^{-\kappa x}}=\frac{1}{D_2M_1}=1 \quad \text{and}
\quad
\lim_{x\to\infty}|u^*_2(x)-1|=0.
\end{equation}
Whence
\begin{equation}\label{zzz-z3-1}
    \lim_{x\to\infty}{\bf u}^*(x)={\bf e}_2.
\end{equation}
In fact ${\bf u}^*(x)$ satisfies the following.

\begin{lemma}\label{lem9}
It holds that
\begin{equation}\label{eq-asym-for-u-star}
     \lim_{x\to\infty}\Big|\frac{u_1^*(x)}{e^{-\kappa x}}-1 \Big|=0 \quad \text{and}\quad \lim_{x\to\infty}\Big| \frac{u_2^*(x)-1}{e^{-\kappa x}}-\frac{(\chi_2\mu_2-rb)}{(1-a)+\frac{r}{M_2}-(d-1)\kappa^2-\frac{\chi_2\mu_2\sqrt{\lambda}}{2}B_{\lambda,\kappa} }\Big|=0,
\end{equation}
where $B_{\lambda,\kappa}$ is as in \eqref{B-lambda-kappa-eq}.
\end{lemma}
\begin{proof}
It is clear that the first limit in \eqref{eq-asym-for-u-star} is established in \eqref{zzz-z3-1-1-1}. So, it remains to show that the second limit in \eqref{zzz-z3-1-1-1} holds. We proceed by contradiction and suppose that there is a sequence $\{x_n\}_{n\ge 1}$ with $x_n\to\infty$ as $n\to\infty$ such that
\begin{equation}\label{j-01}
    \inf_{n\geq 1}\Big| \frac{u_2^*(x_n)-1}{e^{-\kappa x_n}}-\frac{(\chi_2\mu_2-rb)}{(1-a)+\frac{r}{M_2}-(d-1)\kappa^2-\frac{\chi_2\mu_2\sqrt{\lambda}}{2}B_{\lambda,\kappa} }\Big|>0.
\end{equation}
Consider the functions
$$
w_{2}^{n}(x+x_n)=\frac{u_2^*(x_n+x)-1}{e^{-\kappa (x_n+x)}} \quad \text{and}\quad  w_{1}^n(x)=\frac{u_1^*(x_n+x)}{e^{-\kappa(x_n+x)}}, \quad x\in\mathbb{R}
$$
and note that
$$
\|w_{1}^n\|_{\infty}\leq 1 \quad \text{and}\quad \|w_{2}^n\|_{\infty}\leq M_2\tilde{D}_2\quad \forall\ n\ge 1,
$$
since ${\bf u}^*\in\mathcal{E}(\kappa)$. A simple computation shows that $\{(w_{1}^n,w_{2}^n)\}_{n\ge 1}$ satisfy
\begin{align}\label{j-02}
   0=&d(\kappa^2w_2^n-2\kappa w_{2,x}^n+w_{2,xx}^n)+(c_{\kappa}-\chi_2v_{x}(\cdot+x_n;\mu_1u_1^*+\mu_2u_2^*))(w_{2,x}^n-\kappa w_{2}^n)\cr
   & -\left(\frac{rw^n_{2}}{M_2}+(br-\chi_2\mu_1)w^n_1+\chi_2\lambda \tilde{v}(\cdot+x_n;\mu_1w^n_1+\mu_2w^n_2) \right)u_2^{*}(\cdot+x_n)
\end{align}
where the linear and bounded operator $C^{b}_{\rm unif}(\mathbb{R})\ni g\mapsto\tilde{v}(\cdot;g)$ is given by
\begin{equation}\label{j-03}
\tilde{v}(x;g)=\frac{1}{2\sqrt{\lambda}}\int_{\mathbb{R}}e^{-\sqrt{\lambda}|z|-\kappa z}g(z+x)dz \quad \forall\ g\in C^{b}_{\rm unif}(\mathbb{R}).
 \end{equation}
 By a priori-estimates for parabolic equations, we may suppose that $(w_1^n(x),w_2^n(x))\to (w^{\infty}_1(x),w^{\infty}_2(x)) $ as $n\to\infty$ locally uniformly in $C^2_{loc}(\mathbb{{R}})$. Furthermore, we note from \eqref{zzz-z3-1-1-1} that $\lim_{n\to\infty}v(x+x_n;\mu_1u_1^*+\mu_2u_2^*)=\frac{\mu_2}{\lambda}$  and $w^{\infty}_1(x)=1$  for every $x\in\mathbb{R}$. Hence $w^{\infty}_2(\cdot) \in C^{2,b}_{\rm unif}(\mathbb{R})$ and satisfies
 \begin{align}\label{j-04}
   0=&d(\kappa^2w_2^{\infty}-2\kappa w_{2,x}^{\infty}+w_{2,xx}^{\infty})+c_{\kappa}(w_{2,x}^{\infty}-\kappa w_{2}^{\infty}) -\frac{rw_{2}^{\infty}}{M_2}-(br-\chi_2\mu_1)-\chi_2\lambda \tilde{v}(x;\mu_2w_2^{\infty}) \cr
   =&dw_{2,xx}^{\infty}+(c_{\kappa}-2\kappa d)w_{2,x}^{\infty}-\big((1-a)-(d-1)\kappa^2+\frac{r}{M_2})w_2^{\infty}-\chi_2\mu_2\lambda \tilde{v}(x;w_2^{\infty}) -(br-\chi_2\mu_1).
\end{align}
We denote by $W(\cdot;g)$ the solution of the elliptic equation
$$
0=dW_{xx}+(c_{\kappa}-2\kappa d)W_{x}-\big((1-a)-(d-1)\kappa^2+\frac{r}{M_2})W+g(x), \quad x\in\mathbb{R},
$$
for every $g\in C^{b}_{\rm unif}(\mathbb{R})$. Since $(1-a)-(d-1)\kappa^2+\frac{r}{M_2}>0$, it follows from standard theory of elliptic operators, that $W(\cdot;g)\in C^{2,b}_{\rm unif}(\mathbb{R})$ is uniquely determined by $g$. Moreover, the maximum principle implies that
\begin{equation}\label{j-05}
\|W(\cdot;g)\|_{\infty}\leq \frac{1}{(1-a)-(d-1)\kappa^2+\frac{r}{M_2}}\|g\|_{\infty}\quad \forall\ g\in C^{b}_{\rm unif}(\mathbb{R}).
\end{equation}
Now, observe from \eqref{j-04} that
$$
w^{\infty}_{2}(\cdot)=-\chi_2\mu_2\lambda W(\cdot;\tilde{v}(\cdot;w_2^\infty) )-W(\cdot;(br-\chi_2\mu_1)).
$$
Equivalently, we have
\begin{equation}\label{j-05-1}
w^{\infty}_{2}(\cdot)+\chi_2\mu_2\lambda W(\cdot;\tilde{v}(\cdot;w_2^\infty) )=-W(\cdot;(br-\chi_2\mu_1))
\end{equation}
Observe from \eqref{j-03}  that
\begin{equation}\label{j-06}
\|\tilde{v}(\cdot;g)\|_{\infty}\leq \frac{\|g\|_{\infty}}{2\sqrt{\lambda}}\int_{\mathbb{R}}e^{-\sqrt{\lambda}|z|-\kappa z}dz=\frac{B_{\lambda,\kappa}}{2\sqrt{\lambda}}\|g \|_{\infty} \quad\ \forall\ g\in C^{b}_{\rm unif}(\mathbb{R}),
\end{equation}
where $B_{\lambda,\kappa}$ is as in \eqref{B-lambda-kappa-eq}. Hence, by \eqref{j-05} and \eqref{j-06},
$$
\|\chi_2\mu_2\lambda W(\cdot;\tilde{v}(\cdot;g))\|_{\infty}\leq\frac{\chi_2\mu_2\sqrt{\lambda} B_{\lambda,\kappa}}{2\big((1-a)-(d-1)\kappa^2+\frac{r}{M_2} \big) }\|g\|_{\infty} \quad \forall\ g\in C^{b}_{\rm unif}(\mathbb{R}).
$$
We remark from \eqref{xxx-2} that $\frac{\chi_2\mu_2\sqrt{\lambda} B_{\lambda,\kappa}}{2\big((1-a)-(d-1)\kappa^2+\frac{r}{M_2} \big) } <1$. Hence $1$ belongs to the resolvent set of the linear bounded operator $ C^{b}_{\rm unif}(\mathbb{R})\ni g\mapsto -\chi_2\mu_2\lambda W(\cdot;\tilde{v}(\cdot;g))$. As a result, we obtain that the solution of the equation \eqref{j-05-1}, equivalently solution of \eqref{j-04}, is uniquely determined in $C^{2,b}_{\rm unif}(\mathbb{R})$. However, it is easily verified  that the constant function
$$
w(x)=\frac{(\chi_2\mu_2-rb)}{(1-a)+\frac{r}{M_2}-(d-1)\kappa^2-\frac{\chi_2\mu_2\sqrt{\lambda}}{2}B_{\lambda,\kappa} }\quad \forall\ x\in\mathbb{R}
$$
is a solution of \eqref{j-04} in $C^{2,b}_{\rm unif}(\mathbb{R})$. Thus we conclude that
$ w^{\infty}_{2}(x)\equiv \frac{(\chi_2\mu_2-rb)}{(1-a)+\frac{r}{M_2}-(d-1)\kappa^2-\frac{\chi_2\mu_2\sqrt{\lambda}}{2}B_{\lambda,\kappa} }$, in particular we obtain
$$
\frac{(\chi_2\mu_2-rb)}{(1-a)+\frac{r}{M_2}-(d-1)\kappa^2-\frac{\chi_2\mu_2\sqrt{\lambda}}{2}B_{\lambda,\kappa} }=w^{\infty}_2(0)=\lim_{n\to\infty}w_2(x_n),
$$
which contradicts with \eqref{j-01}. Therefore \eqref{eq-asym-for-u-star} must hold.
\end{proof}

Thanks to Lemma \ref{lem9}, to complete the proof of Theorem \ref{main-exist-thm}, it remains to study the asymptotic behavior of $ {\bf u}^*(x)$ and $x\to-\infty$, which we complete in  next two lemmas.

\medskip

\begin{lemma}\label{lem6}
Suppose that  hypothesis ${\bf (H2)}$ holds. Then
$$
\lim_{x\to-\infty}{\bf u}^*(x)={\bf e}_1.
$$
\end{lemma}
\begin{proof}
We prove this result by contradiction. Suppose that there is a sequence $\{x_n\}$ with $x_n\to-\infty$ such that
\begin{equation}\label{zzz-z10}
\inf_{n\geq 1}|{\bf u}^*(x_n)-{\bf e}_1|>0.
\end{equation}
Consider the sequence  ${\bf u}^{*,n}(x)={\bf u}^*(x+x_n)$ for every $x\in\mathbb{R}$ and $n\geq 1$. By a priori estimates for elliptic equations and  Arzela-Ascoli's theorem, without loss of generality, we may suppose that there is some $\tilde{{\bf u}}\in C^2$ such that $ {\bf u}^{*,n} \to \tilde{{\bf u}}$ as $n\to\infty$ locally uniformly in $C^2(\mathbb{R})$. Moreover, $\tilde{{\bf u}} $ also satisfies \eqref{zzz-z2}.  Recalling $m_1^*$ and $x_1$ given by Lemma \ref{lem5} $(ii)$, we deduce that
\begin{equation}\label{zzz-z4}
    m_1^*\leq \tilde{u}_1(x)\leq M_1,\quad \forall\ x\in\mathbb{R}
\end{equation}
since $x_n\to-\infty$ as $n\to\infty$. Therefore by Lemma \ref{lem-01}, we obtain $\tilde{\bf u}(x)\equiv (1,0)$. In particular, $\tilde{\bf u}(0)=e_1$, which contradicts with \eqref{zzz-z10}, since $\tilde{\bf u}(0)=\lim_{n\to\infty}{\bf u}^{*}(x_n)$.
\end{proof}

\begin{lemma}\label{lem7}
Suppose that hypothesis ${\bf (H3)}$ holds. Then
$$
\lim_{x\to-\infty}{\bf u}^*(x)={\bf e}_*.
$$
\end{lemma}
\begin{proof} We proceed also by contradictions. The ideas are similar to that of the of the proof of Lemma \ref{lem6}. Suppose that there is a sequence $\{x_n\}$ with $x_n\to-\infty$ such that
\begin{equation}\label{zzz-z10'}
\inf_{n\geq 1}|{\bf u}^*(x_n)-{\bf e}_*|>0.
\end{equation}
Consider the sequence  ${\bf u}^{*,n}(x)={\bf u}^*(x+x_n)$ for every $x\in\mathbb{R}$ and $n\geq 1$. By a priori estimates for elliptic equations and the Arzela-Ascoli's theorem, without loss of generality, we may suppose that there is some $\tilde{{\bf u}}\in C^2$ such that $ {\bf u}^{*,n} \to \tilde{{\bf u}}$ as $n\to\infty$ locally uniformly in $C^2(\mathbb{R})$. Moreover, $\tilde{{\bf u}} $ also satisfies \eqref{zzz-z2}.  Recalling the positive constants $m_i^*$ and $x_i$, $i=1,2$, and given by Lemma \ref{lem5} $(ii)-(iii)$, we deduce that
\begin{equation}\label{zzz-z4'}
    m_i^*\leq \tilde{u}_i(x)\leq M_i,\quad \forall\ x\in\mathbb{R}, \ i=1,2
\end{equation}
since $x_n\to-\infty$ as $n\to\infty$. Therefore by Lemma \ref{lem-02}, we obtain $\tilde{\bf u}(x)\equiv e_*$. In particular, $\tilde{{\bf u}}(0)={\bf  e}_*$, which contradicts with \eqref{zzz-z10'}, since $\tilde{\bf u}(0)=\lim_{n\to\infty}{\bf u}^{*}(x_n)$.
\end{proof}

Now we complete the proof of Theorem \ref{main-exist-thm}.

\begin{proof}[Proof of Theorem \ref{main-exist-thm}]
For every $c>c_{\kappa^*}$, there is a unique $0<\kappa_c<\kappa^*$ satisfying $c=c_{\kappa_c}=\frac{\kappa_c^2+1-a}{\kappa_{\kappa_c}}$.  This implies that ${\bf u}(t,x)={\bf u}^{*}(x-ct)$ is a traveling solution of \eqref{eq:1} with speed $c$. Moreover, it follows from \eqref{zzz-z3-1} that ${\bf u}^{*}$ connect ${\bf e}_2$ at right end. Since ${\bf u}^{*} \in\mathcal{E}(\kappa)$, then $u^{*}_{1}(x)>0$ by comparison principle for elliptic equations. This in turn implies that $\|u^{*}_2- 1\|_{\infty}>0$. Thus $ {\bf u}^{*}$ is not a trivial solution of \eqref{eq:1}. Assertion $(i)$ and $(ii)$ of the theorem follows from Lemmas \ref{lem6} and \ref{lem7} respectively.
\end{proof}

\section{Proof of Theorem \ref{main-non-exist-thm}}\label{sec-for-t2}

In this section, we present the proof of nonexistence of nontrivial traveling wave solutions of \eqref{eq:1} with speed $c<c_0^*(=2\sqrt{1-a})$ connecting ${\bf e}_2$ at right end.  Our first step toward the proof of the non-existence is to show that, for any nontrivial traveling solution ${\bf u}(x-ct)$ of \eqref{eq:1} connecting  ${\bf u}(\infty)={\bf e}_2$ at the right end,   there holds $u_{1,x}<0$ for $x\gg 1$.

\begin{lemma}\label{lem8}
Let ${\bf u}(t,x)={\bf u}(x-ct)$ be a nontrivial traveling wave solution of \eqref{eq:1} connecting ${\bf e}_2$ at the right end. Then there is $X_0\gg 1$ such that $u_{1,x}(x)\le 0$ for every $x>X_0$.
\end{lemma}
\begin{proof}
We proceed by contraction. Suppose that the statement of the lemma is false. Then, since $u_1(\infty)=0$ and $u_1(x)>0$ for every $x\in\mathbb{R}$, then there is sequence of local minimum points $\{x_n\}_{n\geq 1}$ of $u_1(x)$ satisfying $x_n\to\infty$ as $n\to\infty$. Since $u_2(\infty)=1$, then $\lim_{n\to\infty}u_2(x_n)=1$. From the representation formula
$$
v(x;{\bf u})=\frac{1}{2\sqrt{\lambda}}\int_{\mathbb{R}}e^{-\sqrt{\lambda}|z|}(\mu_1u_1(z+x)+\mu_2u_2(z+x))dz,
$$
it follows from the dominated convergence theorem that
$$
\lim_{n\to\infty}v(x_n;{\bf u})=\frac{\mu_2}{\lambda}.
$$
Hence
$$
\lim_{n\to\infty}\left(1-(1-\chi_1\mu_1)u_1(x_n)-(a-\chi_1\mu_2)u_2(x_n)-\lambda\chi_1v(x_n;{\bf u})\right)=1-a .
$$
Thus there is $n_0\gg 1$ such that
\begin{align}\label{xx-4}
u_{1,xx}(x_n)+(c-\chi_1v_x(x_n;{\bf u}))u_{1,x}(x_n) =&\left((1-\chi_1\mu_1)u_1(x_n)+(a-\chi_1\mu_2)u_2(x_n)+\lambda\chi_1v(x_n;{\bf u})-1\right)u_1(x_n)\cr
<& \frac{-(1-a)u_1(x_n)}{2}<0,\quad \forall\ n\geq n_0.
\end{align}
Since $\{x_n\}$ is a sequence of local minimum points, we have that $u_{1,xx}(x_n)\ge 0$ and $u_{1,x}(x_n)=0$ for every $n\geq 1$, which clearly contradicts with \eqref{xx-4}. Thus the statement of the Lemma must hold.
\end{proof}

Now, we present the proof of Theorem \ref{main-non-exist-thm}.

\begin{proof}[Proof of Theorem \ref{main-non-exist-thm}]
We prove this result by contradiction. Suppose that \eqref{eq:1} has a nontrivial traveling wave solution ${\bf u}(t,x)={\bf u}(x-ct)$ with speed $c<c_0^*$ connecting ${\bf e}_2$ at the right end.  Choose $q>0$ and $0<\varepsilon\ll 1$ satisfying $\max\{c,0\}+\varepsilon<q<2\sqrt{1-a-\varepsilon}$.
 By Lemma \ref{lem8}, there is $X_0\gg 1$ such that $u_{1,x}(x)\le 0$ for every $x>X_0$. Moreover, since ${\bf u}(+\infty)={\bf e}_2$, we deduce that
$$
\lim_{x\to\infty} \left(1-(1-\chi_1\mu_1)u_1(x)-(a-\chi_1\mu_2)u_2(x)-\chi_1\lambda v(x;{\bf u}) \right)=1-a
\quad \text{and}\quad
\lim_{x\to\infty}v_x(x;{\bf u})=0.
$$
Thus, there is  $X_1\gg X_0$, such that
\begin{equation*}
    |\chi_1 v_{x}(x;{\bf u})|<\varepsilon \quad \text{and}\quad 1-(1-\chi_1\mu_1)u_1(x)-(a-\chi_1\mu_2)u_2(x)-\chi_1\lambda v(x;{\bf u})>1-a-\varepsilon \quad \forall\ x\geq X_1.
\end{equation*}
Hence, it holds that the function $\overline{u}(t,x)=u_{1}(x-(c+\varepsilon)t)$ satisfies that
\begin{align}
\overline{u}_{t}=&{u}_{1,xx}-(\varepsilon+\chi_1v_x)u_{1,x}+(1-(1-\chi_1\mu_1)u_1-(a-\chi_1\mu_2u_2)-\lambda\chi_1v(x-(x+\varepsilon)t;{\bf u}))\overline{u}\cr
\geq & u_{1,xx} + (1-a-\varepsilon)\overline{u},\quad\ x\ge X_1+(c+\varepsilon)t,\ t>0.
\end{align}
A simple computation shows that the function
$$
\underline{u}(t,x)= \sigma e^{-\frac{q}{2}(x-x_1-l-qt)}\cos\Big(\frac{\beta}{2}\big(x-x_1-l-qt\big)\Big), \quad x_1+qt\leq x\le l+x_1+qt
$$
where $l=\frac{\pi}{\beta}$, $\beta=\sqrt{4(1-a-\varepsilon)-q^2}$ and $\sigma=e^{-\frac{lq}{2}}\min_{x_1\leq x\leq x_1+L}\overline{u}(0,x)$, satisfies
\begin{equation*}
  \begin{cases}
  \underline{u}_t= \underline{u}_{xx}+(1-a-\varepsilon)\underline{u}, & x_1+qt\leq x\le l+x_1+qt,  \cr
  \underline{u}(t,x)=0 & x=x_1+qt, x=l+x_1+qt. \cr
  \end{cases}
\end{equation*}
Thus, since $q>c+\varepsilon$, then $(c+\varepsilon)t<qt$ for every $t>0$. Moreover the choice of $\sigma$ guarantees that $\underline{u}(0,x)\leq \overline{u}(0,x)$ for every $x_1\leq x\leq x_1+l$ and $\overline{u}(t,x)>0$ for $x\in\{x_1+qt,x_1+l+qt\}$ for every $t>0$. We now infer to the comparison principle for parabolic equations to conclude that
$$
\underline{u}(t,x)\leq \overline{u}(t,x), \quad \forall x_1+qt<x<x_1+l+qt,\ t>0.
$$
In particular,  taking $x=x_1+\frac{l}{2}+qt$, we get
$$
\sigma e^{\frac{ql}{4}}\cos\big(\frac{\pi}{4} \big)= \underline{u}(t,qt+x_1+\frac{l}{2})\leq \overline{u}(t,x_1+qt+\frac{l}{2})=u_1\Big((q-c-\varepsilon)t+x_1+\frac{l}{2}\Big),\quad \forall\ t>0.
$$
Letting $t\to\infty$ yield $ 0<\sigma e^{\frac{ql}{4}}\cos\big(\frac{\pi}{4} \big)\leq u_{1}(\infty)$, which is impossible since $ u_1(\infty)=0$. Therefore, we conclude that the statement of the theorem must hold.
\end{proof}

\medskip

\section{Proof of Theorem \ref{tm-min-wave}}\label{sec-for-min-wave}

In this section, we present the proof of Theorem \ref{tm-min-wave}. To this end, we first recall some results on the spreading speeds and stability for single species chemotaxis model.

\begin{lemma}\cite{Ham1,Sa2}\label{singgle-species-result}
Consider the single species chemotaxis model
\begin{equation}\label{single-species-eq1}
    \begin{cases}
    u_t=\tilde{d}u_{xx}-\chi(uw_x)_x +u(\tilde{a}-\tilde{b}u) & x\in\mathbb{R}, \ t>0\cr
    0=w_{xx}-\tilde{\lambda}w+\mu u & x\in\mathbb{R}, t>0,
    \end{cases}
\end{equation}
where all the parameters are positive, and let $(u(t,x;u_0),w(t,x;u_0))$ denote the unique nonnegative classical of \eqref{single-species-eq1} for every $u_0\in C^b_{\rm unif}(\mathbb{R})$ with $u_0\ge 0$ defined on a maximal interval of existence $[0,t_{\max,u_0})$. Then the following hold.
\begin{itemize}
    \item[(i)] If $\chi\mu<\tilde{b}$, then $t_{\max,u_0}=+\infty$ and  $\|u\|_{\infty}\leq \max\{\|u_0\|,\frac{\tilde{a}}{\tilde{b}-\chi\mu}\}$ for every $t\ge 0$. Moreover, if $\|u_0\|_{\infty}>0$ then
    $$
    \liminf_{t\to\infty}\inf_{|x|\leq (2\sqrt{ad}-\varepsilon)t}u(t,x)>0\quad \forall\ 0<\varepsilon\ll 1.
    $$
    \item[(ii)] If $2\chi\mu<\tilde{b}$ and $\inf_{x\in\mathbb{R}}u_0(x)>0$ then
    $$
    \lim_{t\to\infty}\|u(t,\cdot)-\frac{\tilde{a}}{\tilde{b}}\|_{\infty}=0.
    $$
\end{itemize}
\end{lemma}

Throughout the rest of this section, we assume that {\bf (H5)} holds.
Note that ${\bf (H5)}$ implies ${\bf (H4)}$.
By the definition of the function $F_2(\kappa,\chi_1,\chi_2)$, we have
\begin{align*}
    1-a>  (d-1)_+(1-a)+ F_2(\kappa,\chi_1,\chi_2) +(1-d)\kappa^2\ge F_2(\kappa,\chi_1,\chi_2) \quad \forall \kappa\in (0,\sqrt{1-a}),
\end{align*}
which means that inequality \eqref{xxx-2} also holds for every $ \kappa\in(0,\sqrt{1-a})$.
 Hence $c^*=c_0^*$.

As a result, to complete the proof of Theorem \ref{tm-min-wave}, it remains to show the existence of a non-trivial traveling wave connecting ${\bf e}_2$ at the right end with minimum speed $c_0^*=2\sqrt{1-a}$.

\begin{proof}[Proof of Theorem \ref{tm-min-wave}]

$(i)$ Suppose that hypotheses  ${\bf (H2)}$ and ${\bf (H5)}$ hold and $r>2\chi_2\mu_2$. Chose a decreasing sequence $\{c_n\}_{n\ge 1}$ such that $c_n \to c_0^{*}$ as $n\to\infty$. For every $n\ge 1$, let ${\bf \tilde{u}}^n={\bf \tilde{U}}^{c_n}(x-c_nt)$ be a traveling wave solution of \eqref{eq:1} connecting ${\bf e}_1$ and ${\bf e}_2$ given by Theorem \ref{main-exist-thm}.  Let
\begin{equation}\label{z-z-z-00}
\tilde{x}_n=\min\{x\in\mathbb{R}\ :\ \tilde{U}_1^{c_n}(x)=\frac{1}{2}\}
\end{equation}
and define $ {\bf U}^{c_n}={\bf \tilde{U}}^{c_n}(x+\tilde{x}_n)$ for every $x\in\mathbb{R}$ and $n\ge 1.$ Note that $U_1^{c_n}$ satisfies
\begin{equation}\label{z-z-z-0}
    U_1^{c_n}(x)\begin{cases}
    \ge \frac{1}{2} & \text{if}\ x\le 0\cr
    =\frac{1}{2} & \text{if}\ x= 0
    \end{cases}
\end{equation}
for every $n\ge 1$. Recall that $\|U_i^{c_n}\|_{\infty}<M_i$ for every $n\geq 1$ and $i=1,2$. Hence, by a priori estimates for elliptic equations, if possible after passing to a subsequence, we may suppose that there is some ${\bf U}\in C^{2,b}(\mathbb{R}) $ such that $({\bf U}^{c_n},V(\cdot;{\bf U}^{c_n}))\to({\bf U},V(\cdot;{\bf U}))$ as $n\to\infty$ in $C^{2,b}_{\rm loc}(\mathbb{R})$. Moreover, ${\bf U}(x) $ satisfies
\begin{equation}\label{z-z-z-1}
    \begin{cases}
    0=U_{1,xx}+(c_0^*-\chi_1V_x(\cdot;{\bf U}))U_{1,x} +U_1(1-(1-\chi_1\mu_1)U_1-(a-\chi_1\mu_2)U_2-\chi_1\lambda V(\cdot;{\bf u})) & x\in\mathbb{R}\cr
    0=dU_{2,xx}+(c_0^*-\chi_2V_x(\cdot;{\bf U}))U_{2,x} +U_2(r-(r-\chi_2\mu_2)U_2-(br-\chi_2\mu_1)U_21-\chi_2\lambda V(\cdot;{\bf u})) & x\in\mathbb{R}\cr
    0=V_{xx}-\lambda V+\mu_1U_1+\mu_2U_2 & x\in\mathbb{R}.
    \end{cases}
\end{equation}
 We note that $U_1(\cdot)$ also satisfies properties \eqref{z-z-z-0}. From this point, we complete the proof of the spatial asymptotic behavior of ${\bf U}(x)$ in the following six steps.

\smallskip

 {\bf Step 1}. In this step, we prove that $U_2(x)>0$ for every $x\in\mathbb{R}$.  Suppose not. Then $u(t,x)=U_1(x-c_*t)$ is a solution of the single species chemotaxis model \eqref{single-species-eq1} with $(\tilde{a},\tilde{b},\mu,\tilde{\lambda},\tilde{d})=(1,1,\mu_1,\lambda,1)$. Since $u(0,x)=U_1(x)\geq \frac{1}{2}$ for every $x\le 0$ (by \eqref{z-z-z-0}), then it follows from Lemma \ref{singgle-species-result} (i) that
$$
\liminf_{x\to\infty}U_1(x)=\liminf_{t\to\infty}u((2-\varepsilon)t-c_0^*t)>0, \quad \forall 0<\varepsilon\ll 1.
$$
Hence we conclude that $\inf_{x\in\mathbb{R}}U_1(x)>0$, which yield that $\inf_{t,x\in\mathbb{R}}u(t,x)>0$. It then follows from Lemma \ref{lem-01} that $U_1(x)\equiv 1$. Clearly, this contradicts with \eqref{z-z-z-0} since $U_1(0)=\frac{1}{2}$. Thus we must have that $U_2(x)>0$ for every $x\in\mathbb{R}$.

\smallskip

 {\bf Step 2}. In this step, we prove that $\liminf_{x\to\infty}U_1(x)=0$. If not, since $U_1(\cdot)$ satisfies \eqref{z-z-z-0}, we would have that $\inf_{x\in\mathbb{R}}U_1(x)>0$. And hence since ${\bf (H2)}$ holds, it follows from Lemma \ref{lem-01} that ${\bf U}(x)\equiv {\bf e}_1$, so $U_2\equiv 0$, which contradicts with {\bf Step 1}. Hence $\liminf_{x\in\mathbb{R}}U_1(x)=0$.

\smallskip

{\bf Step 3}. In this step, we prove that $\limsup_{x\to\infty}U_1(x)=0$. Suppose not. According to {\bf Step 2}, there would exist a sequence of minimum points $\{x_n\}_{n\geq 1}$ of the function $U_1 $ satisfying $x_n\to\infty$ and  $U_1(x_n)\to0$ as $n\to\infty$ with $U_{1,x}(x_n)=0$ and $U_{1,xx}(x_n)\leq 0$ for every $n\geq 0$. Hence, we deduce from \eqref{z-z-z-1} that
\begin{align*}
   0 \geq & U_1(x_n)(1-(1-\chi_1\mu_1)U_1(x_n)-(a-\chi_1\mu_2)U_2(x_n)-\chi_1\lambda V(x_n;{\bf u})) \quad \forall\ n\geq 1.
\end{align*}
In particular, we obtain that
\begin{equation*}
    1\leq (1-\chi_1\mu_1)U_1(x_n)+(a-\chi_1\mu_2)U_2(x_n)+\chi_1\lambda V(x_n;{\bf u}) \quad \forall\ n\geq 1,
    \end{equation*}
since $U_1(x_n)>0$ for $n\geq 1$. Letting $n\to\infty$ in the last inequality and using the facts that $ \|U_2\|_{\infty}\leq M_2$ and $\|V(\cdot,{\bf U})\|_{\infty}\leq \frac{1}{\lambda}(\mu_1M_1+\mu_2M_2)$, we obtain that
$$
1\leq (a-\chi_1\mu_2)M_2+\chi_1(\mu_1M_1+\mu_2M_2)=\chi_1\mu_1M_1+aM_2.
$$
This clearly contradicts with hypothesis ${\bf (H2)}$. Thus we must have that $\limsup_{x\to\infty}U_1(x)=0$.

\smallskip

{\bf Step 4.} In this step, we prove that  $\limsup_{x\to\infty}U_2(x)>0$.  If not, then a slight modification of the proof of Theorem \ref{main-non-exist-thm} shows that $c_0^*\geq 2\sqrt{1-aU_2(\infty)}=2$, which is absurd. Hence, we must have $\limsup_{x\to\infty}U_2(x)>0$.

\smallskip

{\bf Step 5}. In this step, we prove that $\liminf_{x\to\infty}U_2(x)>0$. Suppose not. In this case, since $U_1(\infty)=0 $ and $\limsup_{x\to\infty}U_2(x)>0$, we can repeat the arguments used in {\bf Step 3} for the equation satisfies by $U_2(x)$ to end up with the inequality $r\le \chi_2\mu_2M_2$. This clearly contradicts the fact that $r>2\chi_2\mu_2$.

\smallskip

{\bf Step 6.} In this step, we prove that  $\lim_{x\to\infty}U_2(x)=1$. Suppose not. Then there is sequence $\{y_n\}_{n}$ with $y_n\to\infty$ as $n\to\infty$ such that
\begin{equation}\label{z-z-z-3}
\inf_{n\geq 1}|U_2(y_n)-1|>0.
\end{equation}


By a priori estimates for elliptic equations, without loss of generality, we may suppose that there is some ${\bf U}^*\in C^{2,b}_{\rm unif}(\mathbb{R})$ such that ${\bf U}(x+y_n)\to{\bf U}^*(x)$ as $n\to\infty$ in $C^{2,b}_{\rm loc}(\mathbb{R})$. Note by {\bf Step 5} (respectively {\bf Step 3}) that $\inf_{x\in\mathbb{R}}U_2^*(x)>0$ (respectively $U_1^*\equiv0$). Hence, by Lemma \ref{singgle-species-result} $(ii)$, we conclude that $U^*_2(x)\equiv 1$ since $r>2\chi_2\mu_2$. This contradicts with \eqref{z-z-z-3}. Hence $U_2(\infty)=1$.

Finally, we can employ Lemma \ref{lem6} together with the fact that $U_1(x)\geq \frac{1}{2}$ for every $x\leq 0$ to conclude that $\lim_{x\to-\infty}{\bf U}(x)={\bf e}_1$. This completes the proof of $(i)$ .


\smallskip

$(ii)$ Suppose that  ${\bf (H3)}$ and ${\bf H(5)}$ hold. Note that hypothesis ${\bf (H3)}$ implies that $r>2\chi_2\mu_2$. The proof of the minimal wave in this case  follows similar arguments as in $(i)$. So, we shall provide general ideas of the proof. As  in $(i)$, chose a sequence $\{c_n\}_{n\ge 1}$ such that $c_n \to c_0^*+$ as $n\to\infty$. For every $n\ge 1$, let ${\bf \tilde{u}}^n={\bf \tilde{U}}^{c_n}(x-c_nt)$ be a traveling wave solution of \eqref{eq:1} connecting ${\bf e}_1$ and ${\bf e}_2$ given by Theorem \ref{main-exist-thm}.  Next, we let
\begin{equation}\label{z-z-z-000}
\tilde{x}_n^i=\min\{x\in\mathbb{R}\ :\ \tilde{U}_i^{c_n}(x)=\frac{\min\{1-a,1-b\}}{2(1-ab)}\}, \ i=1,2\quad  \text{and} \quad \tilde{x}_n=\min\{\tilde{x}_n^1,\tilde{x}_n^2\},
\end{equation}
and define $ {\bf U}^{c_n}(x)={\bf \tilde{U}}^{c_n}(x+\tilde{x}_n)$ for every $x\in\mathbb{R}$ and $n\ge 1.$ Note that $U_1^{c_n}$ satisfies
\begin{equation}\label{z-z-z-10}
    U_i^{c_n}(x)
    \ge \frac{\min\{1-a,1-b\}}{2(1-ab)} \quad  \forall\  x\le 0, i=1,2, n\ge 1.
\end{equation}
By a priori estimates for elliptic equations, if possible after passing to a subsequence, we may suppose that there is some ${\bf U}\in C^{2,b}(\mathbb{R}) $ such that $({\bf U}^{c_n},V(\cdot;{\bf U}^{c_n}))\to({\bf U},V(\cdot;{\bf U}))$ as $n\to\infty$ in $C^{2,b}_{\rm loc}(\mathbb{R})$. Moreover ${\bf U}(x)$ satisfies \eqref{z-z-z-1}. It is clear from \eqref{z-z-z-10} that $U_i(x)\geq \frac{\min\{1-a,1-b\}}{2(1-ab)}$ for every $x\leq 0$ and $j=1,2$. Hence we may employ Lemma \ref{lem7} to conclude that $\lim_{x\to-\infty}{\bf U}(x)=e_*$. Since $c_0^*=2\sqrt{1-a}<2$, we can proceed as in the proof of {\bf Step 5} to conclude that $\limsup_{x\to\infty}U_2(x)>0$.  Now, we can proceed as in the proof of {\bf Step 6} by using the fact that ${\bf (H3)}$ to conclude that $\liminf_{x\to\infty}U_2(x)>0$. Next, observe from \eqref{z-z-z-10} that there is some $i_0\in\{1,2\}$ such that $U_{i_0}(0)=\frac{\min\{1-a,1-b\}}{2(1-ab)}$. Hence ${\bf U}(x)\not\equiv {\bf e}_*$. Hence, we can proceed as in the proof of {\bf Step 2} using the stability of ${\bf e}_*$ to conclude that $\liminf_{x\to\infty}U_1(x)=0$. This in turn, as is {\bf Step 3} yield that $U_1(\infty)=0$. As result, since $\liminf_{x\to\infty}U_2(x)>0$ and $r>2\chi_2\mu_3$, we may use Lemma \ref{singgle-species-result} $(ii)$ to conclude that $\lim_{x\to\infty}U_2(x)=1$. This completes the proof of $(ii)$.
\end{proof}

\section{Appendix}
In this section we present the proof of Lemmas \ref{lem1} and \ref{lem2}.

\begin{proof}[Proof of Lemma \ref{lem1}]
Observe that for every $x\in\mathbb{R} $
\begin{align*}
\int_{\mathbb{R}}e^{-\sqrt{\lambda}|x-y|}\underline{u}^\kappa_{2}(y)dy=&\int_{\mathbb{R}}e^{-\sqrt{\lambda}|x-y|}\left(1-\tilde{D}_{2}e^{-\kappa y}\right)_+dy\cr
\geq & \int_{\mathbb{R}}e^{-\sqrt{\lambda}|x-y|}\left(1-\tilde{D}_{2}e^{-\kappa y}\right)dy\cr
=& \frac{2}{\sqrt{\lambda}}-\tilde{D}_{2}\int_{\mathbb{R}}e^{-\sqrt{\lambda}y-\kappa y}dy=\frac{2}{\sqrt{\lambda}}-\tilde{D}_{2}B_{\lambda,\kappa}e^{-\kappa x}.
\end{align*}
It is clear that $\int_{\mathbb{R}}e^{-\sqrt{\lambda}|x-y|}\underline{u}^\kappa_{2}(y)dy> 0 $ for all $x\in\mathbb{R}$. Hence
\begin{equation}\label{z-1}
\int_{\mathbb{R}}e^{-\sqrt{\lambda}|x-y|}\underline{u}^\kappa_{2}(y)dy\geq \left(\frac{2}{\sqrt{\lambda}}-\tilde{D}_{2}B_{\lambda,\kappa}e^{-\kappa x} \right)_+.
\end{equation}
Similarly for every $x\in\mathbb{R}$
\begin{align*}
\int_{\mathbb{R}}e^{-\sqrt{\lambda}|x-y|}\underline{u}^\kappa_{1}(y)dy\ge&M_1D_2\int_{\mathbb{R}}e^{-\sqrt{\lambda}|x-y|}\left(1-D_{1}e^{-\varepsilon_1 y}\right)e^{-\kappa y}dy=M_1D_2e^{-\kappa x}\left(B_{\lambda,\kappa}-D_1B_{\lambda,\kappa+\varepsilon_1} e^{-\varepsilon_1 x}\right).
\end{align*}
Hence, since $\int_{\mathbb{R}}e^{-\sqrt{\lambda}|x-y|}\underline{u}^\kappa_{1}(y)dy> 0$ for every $x\in\mathbb{R}$, we deduce that
$$
\int_{\mathbb{R}}e^{-\sqrt{\lambda}|x-y|}\underline{u}^\kappa_{1}(y)dy\geq M_1D_2e^{-\kappa x}\left(B_{\lambda,\kappa}-D_1B_{\lambda,\kappa+\varepsilon_1} e^{-\varepsilon_1 x}\right)_+,
$$
which together with \eqref{z-1} yields \eqref{xx-2} since  $$v(x;{\bf u})\geq \frac{\mu_1}{2\sqrt{\lambda}}\int_{\mathbb{R}}e^{-\sqrt{\lambda}|x-y|}\underline{u}^\kappa_{1}(y)dy +\frac{\mu_2}{2\sqrt{\lambda}}\int_{\mathbb{R}}e^{-\sqrt{\lambda}|x-y|}\underline{u}^\kappa_{2}(y)dy. $$
\end{proof}

\begin{proof}[Proof of Lemma \ref{lem2}]. For every $x\in\mathbb{R}$, observe from \eqref{eq2} that
\begin{equation}\label{zz-z1}
\frac{d}{dx}v(x;{\bf u})=\frac{1}{2}\int_{\mathbb{R}}\text{sign}(y-x)e^{-\sqrt{\lambda}|y-x|}(\mu_1u_1(y)+\mu_2u_2(y))dy,\quad \forall\ x\in\mathbb{R}.
\end{equation}
Hence, since $0\leq u_1(x)\leq \overline{u}^\kappa_1(x)=\min\{M_1,M_1D_2e^{-\kappa x}\}$, we obtain
\begin{equation}\label{zz-z2}
\left|\int_{\mathbb{R}}\text{sign}(y-x)e^{-\sqrt{\lambda}|y-x|}u_1(y)dy\right|\leq \int_{\mathbb{R}}e^{-\sqrt{\lambda}|y-x|}\overline{u}^\kappa_1(y)dy\leq \min\left\{\frac{2M_1}{\sqrt{\lambda}},M_1D_2B_{\lambda,\kappa}e^{\kappa x}\right\}, \quad \forall\ x\in\mathbb{R}.
\end{equation}
On the other hand using the fact that $\int_\mathbb{R}\text{sign}(z)e^{-\sqrt{\lambda}z}dz$=0, and that $|u_2(x)-1|\le M_2\tilde{D}_{2}e^{-\kappa x}$, we obtain
$$
\left|\int_{\mathbb{R}}\text{sign}(y-x)e^{-\sqrt{\lambda}|y-x|}u_2(y)dy\right|=\left| \int_{\mathbb{R}}\text{sign}(y-x)e^{-\sqrt{\lambda}|y-x|}(\overline{u}^\kappa_1(y)-1)dy\right|\leq M_{2}\tilde{D}_{2}B_{\lambda,\kappa}e^{\kappa x}, \quad \forall\ x\in\mathbb{R},
$$
which combined with the fact that
$$
\left|\int_{\mathbb{R}}\text{sign}(y-x)e^{-\sqrt{\lambda}|y-x|}u_2(y)dy\right|\leq \int_{\mathbb{R}}e^{-\sqrt{\lambda}|y-x|}\|u_2\|_{\infty}dy\leq \frac{2M_2}{\sqrt{\lambda}}, \quad \forall\ x\in\mathbb{R}
$$
yields
\begin{equation}\label{zz-z3}
\left|\int_{\mathbb{R}}\text{sign}(y-x)e^{-\sqrt{\lambda}|y-x|}u_2(y)dy\right|\leq \min\left\{\frac{2M_2}{\sqrt{\lambda}}, M_2\tilde{D}_{2}B_{\lambda,\kappa}e^{\kappa x} \right\}\quad \forall\ x\in\mathbb{R}.
\end{equation}
The statement of Lemma \ref{lem2} follows from \eqref{zz-z1}-\eqref{zz-z3}.
\end{proof}


\begin{thebibliography}{}

\bibitem{Lankeit1} T. Black, J. Lankeit, M. Mizukami, On the weakly competitive case in a two-species
chemotaxis model, {\it  IMA Journal of Applied Mathematics}, {\bf 81} (2016), 860–876.

\bibitem{Fan1} J. Fang, X-Q. Zhao, Traveling waves for monotone semiflows with weak compactness, {\it SIAM J. Math. Anal.},  {\bf 46}, No. 6 (2014),  3678–3704

\bibitem{Gar1} R. A. Gardner, Existence and stability of travelling wave solutions of competition models: a degree theoretic approach,  {\it J. Differential  Equations}, {\bf 44} (1982), 343-364.

\bibitem{Gar2} R. A. Gardner and C. K. R. T. Jones, Stability of travelling wave solutions of diffusive predator-prey systems, {\it Trans. Amer. Math. Soc.}, {\bf 327} (1991), 465-524.

\bibitem{Huang1} W. Huang, M. Han, Non-linear determinacy of minimum wave speed fora Lotka–Volterra competition model, {\it Journal of Differential Equations}, {\bf 251} (2011), 1549-1561.

\bibitem{Huang2} W. Huang, Problem on minimum wave speed for a Lotka–Volterra reaction–diffusion competition model, {\it J. Dynam. Differential Equations}, {\bf 22} (2010), 285–297.

\bibitem{Li2} J-S. Guo,  X. Liang,  The Minimal Speed of Traveling Fronts for the Lotka–Volterra Competition System, {\it Journal of Dynamics and Differential Equations}, {\bf 23} (2011), 353–363.

\bibitem{Trudinger} D. Gilbert and S. N. Trudinger, {Elliptic partial differential equations of second order}, Springer-Verlag, Berlin, 2001.

\bibitem{Ham1} F. Hamel, C. Henderson, Propagation in a Fisher-KPP equation with non-local advection, {\it Journal of Functional Analysis}, {\bf 278} (2020), 108426, 53 pp.

\bibitem{Hartman} Phillip Hartman, Ordinary differential equations, Second Edition, {\it Boston, Basel; Stuttgart : Birkh\"auser, 1982}

\bibitem{Hos1} Y. Hosono, Singular perturbation analysis of travelling waves of dffusive Lotka-Volterra competition models, Numerical and applied mathematics, Part II (Paris, 1988), 687-692, IMACS Ann. Comput. Appl. Math., 1. 2, Baltzer, Basel, 1989

\bibitem{Issa1} T. B. Issa and  R. B. Salako, Asymptotic dynamics in a two-species chemotaxis model with non-local terms, {\it Discrete and Continuous Dynamical Systems Series B}, {\bf 22}  (2017),  3839 - 3874.

\bibitem{Kan1} Y. Kan-on, Parameter dependence of propagation speed of travelling waves for competition-diffusion equations,
{\it SIAM J. Math. Anal.}, {\bf 26} (1995), 340-363.

\bibitem{Kan2} Y.Kan-on, Existence of standing waves for competition-diffusion equations,
{\it Japan J. Indust. Appl. Math.}, {\bf 13} (1996), 117-133.

\bibitem{Kan3} Y. Kan-on, Fisher wave fronts for the Lotka-Volterra competition model with diffusion,
{\it Nonlinear Anal.}, {\bf 28} (1997), 145-164.



\bibitem{Le1} Mark A. Lewis, Bingtuan Li, and Hans F. Weinberger, Spreading speed and linear determinacy for two-species competition models,  {\it J. Math. Biol.}, {\bf 45} (2002), no.3, 219-233.

\bibitem{Li1} Bingtuan Li, Hans F. Weinberger, and Mark A. Lewis, Spreading speeds as slowest wave speeds for cooperative systems, {\it Math. Biosci.}, {\bf 196} (2005), no. 1, 82-98.


\bibitem{Miz1} M. Mizukami, T. Yokota1 Global existence and asymptotic stability of solutions to a two-species chemotaxis system with any chemical diffusion, {\it J. Differential Eq.}, {\bf  261} (2016), 2650-2669.

\bibitem{Neg1} M. Negreanu and J.I. Tello, On a competitive system under chemotaxis effects with nonlocal terms, {\it Nonlinearity}, {\bf 26}  (2013), 1083-1103.

\bibitem{Neg2} M. Negreanu and J.I. Tello, Asymptotic stability of a two species chemotaxis system with non-diffusive chemoattractant, {\it J. Differential Eq.}, {\bf 258} (2015), 1592-1617.

\bibitem{Neg3} M. Negreanu and J. I. Tello, On a two species chemotaxis model with slow chemical diffusion, {\it SIAM J. Math. Anal.}, {\bf 46} (2014), no. 6, 3761-3781.


\bibitem{Nadin1} G. Nadin, B. Perthame, L. Ryzhik,   Traveling waves for the Keller–Segel system with Fisher birth terms, {\it Interfaces and Free Boundaries}, {\bf 10} (2008), 517–538.


\bibitem{Sa1} R. B. Salako, W, Shen, Traveling wave solutions for fully parabolic Keller-Segel chemotaxis systems with a logistic source, {\it Electron. J. Differential Equations}, {\bf 2020} (2020), No. 53, pp. 1-18.

\bibitem{Sa2} R. B. Salako, W. Shen, X. Shuwen, Can chemotaxis speed up or slow down the
spatial spreading in parabolic-elliptic Keller-Segel systems with logistic
source? {\it Journal of Mathematical Biology}, {\bf 79} (2019), no. 4, 1455-1490.

\bibitem{Sa3} R. B. Salako, W. Shen,  Existence of Traveling wave solution of
parabolic-parabolic chemotaxis systems, {\it  Nonlinear Analysis: Real World
Applications}, {\bf 42} (2018),  93-119.

\bibitem{Sa4} R. B. Salako, W. Shen, Spreading Speeds and Traveling waves of a
parabolic-elliptic chemotaxis system with logistic source on $\mathbb{R}^N$, {\it Discrete and Continuous Dynamical Systems - Series A}, {\bf 37} (2017),
6189-6225.

\bibitem{Tello1} C. Stinner, J.I. Tello, and W. Winkler, Competive exclusion in a two-species chemotaxis, {\it J.Math.Biol.},
{\bf 68} (2014), 1607-1626.



\bibitem{fife} M. M. Tang and P. C. Fife, Propagating fronts for competing species equations with diffusion,
 {\it Arch. Rational Mech. Anal.}, {\bf 73} (1980), 69-77.

\bibitem{Tello2} J.I. Tello and M. Winkler, Stabilization in two-species chemotaxis with a logistic source, {\it Nonlinearity}, {\bf 25} (2012), 1413-1425.




\end{thebibliography}
\end{document}